\documentclass[11pt]{amsart}
\usepackage{amsmath}
\usepackage{amsmath, amsthm, amscd, amsfonts, amssymb,color, tikz, enumerate}
\usepackage{amsfonts}
\usepackage{amsmath,amsthm,amscd,amssymb,mathrsfs,setspace}
\usepackage{latexsym,epsf,epsfig}
\usepackage{bm}
\usepackage{color}
\usepackage[hmargin=2.5cm,vmargin=3cm]{geometry}
\usepackage[backref]{hyperref}
\usepackage{cite}
\usepackage{mathtools}
\usepackage{dirtytalk}
\usepackage{cancel}
\usepackage[shortlabels]{enumitem}
\usepackage [english]{babel}
\usepackage [autostyle, english = american]{csquotes}
\usepackage{mathtools}
\MakeOuterQuote{"}

\newtheorem{theorem}{Theorem}[section]
\newtheorem{lemma}[theorem]{Lemma}
\newtheorem{proposition}[theorem]{Proposition}

\newtheorem{corollary}[theorem]{Corollary}

\newtheorem{definition}{Definition}

\theoremstyle{remark}
\newtheorem{remark}{Remark}[section]
\numberwithin{equation}{section}

\DeclarePairedDelimiter\floor{\lfloor}{\rfloor}

\theoremstyle{plain}

\newcommand{\comments}[1]{}

\newcommand{\R}{\mathbb R}

\newcommand{\D}{\displaystyle }
\newcommand{\dt}{{\D\frac{d}{dt}}}
\newcommand{\ra}{\rightarrow}
\newcommand{\lra}{\longrightarrow}
\newcommand{\be}{\begin{equation}}
\newcommand{\ee}{\end{equation}}
\newcommand{\bes}{\begin{equation*}}
\newcommand{\ees}{\end{equation*}}

\newcommand{\bea}{\begin{eqnarray}}
\newcommand{\eea}{\end{eqnarray}}
\newcommand{\beas}{\begin{eqnarray*}}
\newcommand{\eeas}{\end{eqnarray*}}

\def\cal#1{\mathcal{#1}}

\def\<{\langle} 
\def\>{\rangle}

\newcommand{\h}{\mathbb H}

\title[Determining Map and Data Assimilation for 3D Fluids]{Determining Map, Data Assimilation and an Observable Regularity Criterion for the Three-Dimensional Boussinesq System}
\author{Abhishek Balakrishna}

\author{Animikh Biswas}

\makeatletter
\newcommand*{\rom}[1]{\expandafter\@slowromancap\romannumeral #1@}
\makeatother

\begin{document}
\begin{abstract}
In this paper, we provide conditions, \emph{based solely on the observed velocity data},  for the global well-posedness, regularity and convergence of the Azouni-Olson-Titi data assimilation algorithm (AOT algorithm) for  a Leray-Hopf weak solutions  of the three dimensional Boussinesq system.  This condition also guarantees the construction of the {\it determining map}. The aforementioned conditions on the (finite-dimensional) velocity observations, which in this case comprise either of a finite-dimensional \emph{modal} projection or finitely many \emph{volume element observations}, are automatically satisfied for solutions that are globally regular and are uniformly bounded in the $H^1$-norm. However, neither regularity nor uniqueness is {\it a priori} assumed on the solutions. To the best of our knowledge, this is the first such rigorous analysis of the AOT  data assimilation algorithm for the  three-dimensional Boussinesq system. As a corollary, we obtain 
that the condition that we imposed in order to obtain well-posedness and tracking property for the AOT system 
is in fact {\it a new observable regularity criterion on the weak global attractor.} The proof of this fact proceeds through the construction of the determining map.
\end{abstract}
		\maketitle
		\section{Introduction}
		For a given dynamical system, which is believed to accurately describe some aspect(s) of an underlying physical reality, 
the problem of forecasting is often hindered by  inadequate knowledge of the initial state and/or model parameters describing the system. However, in many cases, this is compensated by the fact that one has access to data from (possibly noisy) measurements of the system, collected  on a {\it much coarser spatial grid than the desired resolution of the forecast.} The objective of data assimilation and signal synchronization is to use this coarse scale observational measurements to fine tune  our knowledge of the state and/or model to improve the accuracy of the forecasts \cite{Daley1991, Kalnay2003}.

Due to its ubiquity in scientific applications, data assimilation has been the subject of a very large body of work.  Classically, these techniques are based on linear quadratic estimation, also known as the Kalman Filter.  The Kalman Filter has the drawback of assuming that the underlying system and any corresponding observation models are linear.  It also assumes that measurement noise is Gaussian distributed.  This has been mitigated by practitioners via modifications, such as the Ensemble Kalman Filter, Extended Kalman Filter and the Unscented Kalman Filter and consequently, there has been a recent surge of interest in developing a rigorous mathematical framework for these approaches; see, for instance, \cite{asch,  HMbook2012,  Kalnay2003, KLS, LSZbook2015, ReichCotterbook2015} and the references therein. These works provide a Bayesian and variational framework for the problem, with emphasis on analyzing variational and Kalman filter based methods. It should be noted however that the problems of stability, accuracy and {\it catastrophic filter divergence}, particularly for {\it infinite dimensional chaotic dynamical systems governed by PDE's,} continue to pose serious challenges to rigorous analysis, and are far from being resolved \cite{HMbook2012, BLSZ2013, BLLMCSS2013, tmk2016-1, tmk2016-2}.

An alternative approach to data assimilation, henceforth referred to as the \emph{AOT algorithm,}
has recently been  proposed in \cite{AzouaniTiti2014, AOT}, 
which employs a feedback control paradigm
via a \emph{Newtonian relaxation scheme} (\emph{nudging}). 
This was in turn predicated on the notion of finite determining functionals (modes, nodes, volume elements)
for dissipative systems, the rigorous  existence of which was first
established in  \cite{FoiasProdi1967, FoiasTemam1984, FoTi91, JT92}. Assuming that the observations are generated from a continuous dynamical system given by
\[
\dt u = F(u), u(0)=u_0,
\]
the AOT algorithm entails solving an associated system
\begin{equation} \label{aotsystem}
\dt w =F(w)-\mu (I_hw-I_hu), w(0)=w_0\ (\mbox{arbitrary}),
\end{equation}
where $I_h$ is a finite rank linear operator acting on the phase space, called \emph{interpolant operator}, constructed \emph{solely from observations} on $u$ (e.g. low (Fourier) modes of $u$ or values of $u$ measured in a coarse spatial grid). Here $h$ refers to the size of the spatial grid or, in case of the \emph{modal interpolant}, the reciprocal of $h$ stands for the number of observed modes. Moreover,  $\mu>0$ is the \emph{relaxation/nudging parameter} an appropriate choice of which needs to be made for the algorithm to work. It can then be established that the AOT system
\eqref{aotsystem} is well-posed and its solution \emph{tracks} the solution of the original system asymptotically, i.e. $\|w-u\| \lra 0$ as $t \ra \infty$ in a suitable norm.

Although initially introduced in the context of the two-dimensional Navier-Stokes equations, this was later generalized to include various other models and convergence in stronger norms (e.g. the analytic Gevrey class)  \cite{ANT,BiswasMartinez2017,
FGMW, FarhatJollyTiti2015, FLT, FLT1, FLT2, MarkowichTitiTrabelsi2016, pei}, as well as to more general situations such as discrete in time and  error-contaminated measurements and to statistical solutions  \cite{BessaihOlsonTiti2015, IHP, FoiasMondainiTiti2016}. This method has  been shown to perform remarkably well in numerical simulations \cite{AltafTitiKnioZhaoMcCabeHoteit2015, FJJT, GeshoOlsonTiti2015, HOT, HJ2018, larios2019} and has recently been successfully implemented for the first time for efficient dynamical downscaling of a global atmospheric reanalysis \cite{desam2019}. Recent applications include its implementation in \emph{reduced order modeling (ROM)} of turbulent flows to mitigate inaccuracies in ROM \cite{Rebholz2019}, and in
inferring flow parameters and turbulence configurations \cite{DiLeoni1, CHL}.

In this paper, we consider the well-posedness, stability and convergence/tracking property of solutions of the AOT system for Leray-Hopf weak solutions of the three dimensional Boussinesq system, under adequate conditions on the {\it observations of the velocity field, and without assuming regularity}. This extends our previous work on the Navier-Stokes equations in \cite{Nav3d}. To the best of our knowledge, this is the first such rigorous result for the 3D Boussinesq equations.
In all the cases mentioned before except \cite{Nav3d} where rigorous analysis is available, including the Navier-Stokes-$\alpha$ models \cite{AB, ANT, FLT2}, one crucially uses the fact that these models are well-posed and regular, i.e. uniform-in-time bound in a higher Sobolev norm (e.g. the $\h^1$-norm) is available. {\it These bounds are  used in providing an upper bound on the spatial resolution $h$} of the observations (or lower bound on the number of observed low modes)
necessary for the algorithm to be well-posed, stable and convergent. Additionally, the value of the nudging parameter $\mu$ guaranteeing convergence/tracking property also explicitly depends on this uniform bound. By contrast, our condition is formulated purely in terms of the observed, finite dimensional part of the data. For $(u,\theta)$, the solution to the Boussinesq system given in \eqref{main1}-\eqref{main3}, the observations on the velocity component $u$ are used to define the quantity
\begin{equation}\label{tsk}
     \begin{split}
        M_{h,u}^2=32\sup_{0\leq t\leq T}\begin{dcases}
        \displaystyle\|P_N(u)\|^2\sim \sum_{|k|\leq N}|\lambda_{k}|^2|\hat{u}(k)|^2,~N\sim\frac{1}{h}&(\text{Modal})\\[10pt]
         \displaystyle Ch \sum_{\alpha}|\bar{u}_{\alpha}|^2,
          ~\bar{u}_\alpha = \frac{1}{|Q_\alpha|}\int_{Q_\alpha} u &(Volume)
        \end{dcases}
        \end{split}
\end{equation}
where, $\lambda_{k}$ are eigenvalues of the Stokes operator such that $\lambda_{1}\leq\lambda_{2}\leq\dots$ with $P_N$ being the associated spectral projections. Also, in the above definition, $\{Q_\alpha\}$ 
denotes partition of the domain into finitely many cubes of side length $h$. $M_h$ so defined is a finite quantity (as shown in Remark \ref{finMh}). Then, if for $h_0>0$ (depending on fluid viscosity, thermal diffusivity and the size of the domain) there exists $0< h\leq h_0$ such that 
\begin{equation}\label{dai}
    \frac{16cM_{h,u}^4}{\nu^3}\leq\frac{\nu}{4ch^2},\quad\text{($\nu$ is the fluid viscosity)}
\end{equation}
then a choice of $\mu$  exists such that the data assimilated solution $(w,\eta)$ is \emph{regular} (given in Definition \ref{definesol}) and converges to the actual solution $u$. 
\subsection{A New Observanle Regularity Criterion}
Condition \eqref{dai} always holds when the velocity component $u$ is \emph{regular} (i.e. $\displaystyle \sup_{t \in [0,T)}\|u\|_{H^1(\Omega)}<\infty$). So one wonders about the relationship between \eqref{dai} and the regularity of the solution. This is explored in Section \ref{attractor} (and section \ref{weakenn}) to obtain a {\it new observable regularity criteria} for the 3D Navier-Stokes on the weak attractor; similar result can be obtained for the three dimensional Boussinesq system as well. To describe our result,  let $h_0>0$ be defined as
         \[
             h_0^{-2}=\max\left\{\frac{1}{4c\lambda_1},\frac{32c|f|^4}{\nu^8\lambda_1^2}\right\}\quad \text{($f$ is the body force)}.
         \]
         Let $u(t), t \in \R$ be a Leray-Hopf weak solution of the 3D Navier-Stokes  on the weak global attractor $\cal A$. Let $M_h$ be defined as in \eqref{tsk}, except that the supremum is taken on the interval $(-\infty,T]$.
         Assume there exists $0<h\leq h_0$ for which $$\frac{2cM_h^4}{\nu^3}\leq \frac{\nu}{16ch^2}.$$ Then $u(t)$ is regular on $(-\infty,T]$. We refer to this as {\it an observable regularity criterion} because it is purely based on the finitely many observations (modes or volume elements).

	\section{3-D Boussinesq Equation}
	\subsection{Model}
	The Bénard convection problem is a model of the Boussinesq convection system of an incompressible fluid layer, confined between two solid walls, which is heated from below in such a way that the lower wall maintains a temperature $T_0$, while the upper one maintains a temperature $T_1 < T_0$. In this case, after some change of variables and proper scaling (by normalizing the distance between the walls and the temperature difference), the three-dimensional Boussinesq equations that govern the perturbation of the velocity($u$) and temperature about the pure conduction steady state($\theta$) are
	\begin{align}
	&\frac{du}{dt}+\nu \Delta u+(u\cdot\nabla)u=\theta \bm{e}_3 \label{B1}\\ 
	&\frac{d\theta}{dt}+\kappa \Delta \theta+(u\cdot\nabla)\theta- u\cdot\bm{e}_3=0 \label{B2}\\ 
	&\nabla \cdot u=0\\
	& u(0,x)=u_0,~~\theta(0,x)=\theta_0 \label{B3}
	\end{align} 
		For boundary conditions, in the $x_3$ direction we have $$u,\theta=0, \text{ at } x_3=0 \text{ and } x_3=1 $$ and in the $x_1$ and $x_2$ directions, for simplicity, we have the periodic boundary condition $$u,\theta \text{ are periodic, of period $L$ in the $x_1$ and $x_2$ directions.}$$ Here, $x=(x_1,x_2,x_3)$ is a point in the domain, $u(t;x) = (u_1(t, x), u_2(t, x),u_3(t,x))$ is the fluid velocity and $\theta = \theta (t, x)$ is the scaled fluctuation of the temperature around the pure scaled conduction steady-state background temperature profile $1 - x_3$. It is given by $\theta = T - (1 - x_3)$, where $T = T (t, x)$ is the scaled temperature of the fluid inside the domain  $\Omega$. $\kappa$ and $\nu$ are the thermal diffusivity and kinematic viscosity, respectively.\\
    
	\subsection{Notation}	
	We briefly introduce the functional setting for \eqref{B1}-\eqref{B3}. For $\alpha>0$, $H^{\alpha}(\Omega)$ is the usual Sobolev space. We denote the inner product and norm of $L^2(\Omega)$ by $(\cdot,\cdot)$ and $|\cdot|$ respectively and the inner product and norm of $H^1(\Omega)$ by $((\cdot,\cdot))$ and $||\cdot||$ respectively. We define $\mathcal{F}$ to be the set of $C^\infty(\Omega)$ functions defined in $\Omega$, which are trigonometric polynomials in $x_1$ and $x_2$ with period $L$, and compactly supported in the $x_3$-direction. We denote the space of vector valued functions on $\Omega$ that incorporates the divergence free condition by
	$\mathcal{V}=\left\{\phi\in \mathcal{F}\times\mathcal{F}\times\mathcal{F}|\nabla . \phi=0\right\}.$ $H_0$ and $H_1$ are closures of $\mathcal{V}$ and $\mathcal{F}$  in $L^2(\Omega)$ respectively and $V_0$ and $V_1$ are closures of $\mathcal{V}$ and $\mathcal{F}$ in $H^1(\Omega)$ respectively.\\
	$H_0$ and $H_1$ are endowed with the inner products
	$$(u,v)_0=\sum_{i=1}^3\int_{\Omega}u_i(x)v_i(x)dx$$ and $$(\phi,\psi)_1=\int_{\Omega}\phi(x)\psi(x)dx$$ respectively, and the norms $|u|_0=(u,u)_0^{1/2}$ and $|\phi|_1=(\phi,\phi)_1^{1/2}$ respectively. $V_0$ and $V_1$ are endowed with the inner products
	$$((u,v))_0=\sum_{i,j=1}^3\int_{\Omega}\partial_ju_i(x)\partial_jv_i(x)dx,$$
	$$((\phi,\psi))_1=\sum_{j=1}^3\int_{\Omega}\partial_j\phi(x)\partial_j\psi(x)dx$$ respectively, and the associated norms $\|u\|_0=((u,u))_0^{1/2}$ and $\|\phi\|_1=((\phi,\phi))_1^{1/2}$ respectively. We also denote by $P_{\sigma}$ the Leray-Hopf orthogonal projection operator from $L^2(\Omega)$ to $H_0$.\\
	
	\noindent
	Let $D(A_0)=V_0\cap \left(H^2(\Omega)\right)^3$, $D(A_1)=V_1\cap H^2(\Omega)$ and $A_i : D(A_i)\to H_i$ be the unbounded linear operator defined by
	\begin{equation}\label{A}
	(A_iu,v)_i=((u,v))_i,~~\text{for }i=0,1.
	\end{equation}
	We recall that $A_i,$ for $i=0,1,$ is a positive self adjoint operator with a compact inverse. Moreover, there exists a complete orthonormal set of
	eigenfunctions $\phi_{j,i}\in H_i $, such that $A_i\phi_{j,i} = \lambda_{j,i}\phi_{j,i}$, where $0<\lambda_{1,i} \leq \lambda_{2,i} \leq \lambda_{3,i} \leq \dots $ are the
	eigenvalues of $A_i$ repeated according to multiplicity.\\
	
	\noindent
	For $i=0,1,$ we denote by $H_n^i$ the space spanned by the first $n$ eigenvectors of $A_i$ and the
	orthogonal projection from $H_i$ onto $H_n^i$ is denoted by $P_n^i$. We also have the
	Poincare inequality
	\begin{equation}\label{poincare}
	\lambda_1^{1/2}|v|_i\leq\|v\|_i, v\in V_i.
	\end{equation}
	where $\lambda_1=\min\left\{\lambda_{1,0},\lambda_{1,1}\right\}.$\\
	
	\noindent
	Let $V_i'$ be the dual of $V_i$ for $i=0,1$. We define the bilinear term $B_0:V_0\times V_0\to V_0'$ by
	\[
	\langle B_0(u,v),w\rangle_{V_0',V_0}=(((u\cdot\nabla) v),w)_0
	\]
	and $B_1:V_0\times V_1\to V_1'$ by 
	\[
	\langle B_1(u,v),w\rangle_{V_1',V_1}=(((u\cdot\nabla) v),w)_1.
	\]
		
	\noindent
	The bilinear term $B_i$, for $i=0,1,$ satisfies the orthogonality property
	\begin{equation}\label{orthoganal}
	B_i(u,w,w) = 0~ \forall~ u\in V_0,~w\in V_i.
	\end{equation}
	We define the solution space $H=H_0\times H_1$ equipped with the inner product $$\langle s_1,s_2\rangle=(u_1,u_2)_0+(\theta_1,\theta_2)_1,$$ where $s_1=(u_1,\theta_1)$ and $s_2=(u_2,\theta_2).$\\
	We recall some well-known bounds on the bilinear term for velocity in the 3D case.
	\begin{proposition}
		If $u,v\in V_0$ and $w\in H_0$, then
		\begin{equation}\label{nolinest1}
		|(B_0(u,v),w)_0|\leq c\|u\|_{L^6}\|\nabla v\|_{L^3}\|w\|_{L^2}\leq c\|u\|\|v\||Av|^{1/2}|w|
		\end{equation}
		Moreover if $u,v,w\in V_0$, then
		\begin{equation}\label{nolinest2}
		|(B_0(u,v),w)_0|\leq c\|u\|_{L^4}\|\nabla v\|_{L^2}\|w\|_{L^4}\leq c|u|^{1/4}\|u\|^{3/4}\|v\||w|^{1/4}\|w\|^{3/4}
		\end{equation}
	\end{proposition}
	\noindent
We also recall the Ladyzhenskaya's inequality for three dimensions :
\begin{equation}\label{lady}
\|w\|_{L^4}\leq C|w|_0^{\frac{1}{4}}\|w\|_0^{\frac{3}{4}}
\end{equation}
\begin{remark}
	From \cite{fmtr}, we know that the weak solution is uniformly bounded in time in $H=H_0\times H_1$. Hence, for $u$ and $\theta$ as in \eqref{B1}-\eqref{B2}, there exists $M_0(u_0,\theta_0)$, $M_1(u_0,\theta_0)\in \mathbb{R}$ such that
	\begin{equation}\label{unibound}
	|u|_0\leq M_0, \text{ and } ~|\theta|_1\leq M_1.
	\end{equation}
\end{remark}
	
\noindent
We denote by $P_{\sigma}$ the Leray-Hopf orthogonal projection operator from $L^2(\Omega)$ to $H_0$.
With above notation, by applying $P_{\sigma}$ to \eqref{B1}, we may express the 3-D Boussinesqu equation in the following functional form:
	\begin{align}
	&\frac{du}{dt}+\nu A_0(u)+B_0(u,u)=P_{\sigma}(\theta \bm{e}_3) \label{main1}\\ 
	&\frac{d\theta}{dt}+\kappa A_1(\theta)+B_1(u,\theta)- u\cdot\bm{e}_3=0 \label{main2}\\ 
	& u(0,x)=u_0,~~\theta(0,x)=\theta_0 \label{main3}.
	\end{align} 
		\begin{definition}\label{defweak0}
	    $(u,\theta)$ is said to be a weak solution to \eqref{main1}-\eqref{main3} if for all $T>0$,\\
	    \begin{itemize}
	        \item $u\in L^{\infty}(0,T;H_0)\cap L^{2}(0,T;V_0)$\\
	        \item $\theta\in L^{\infty}(0,T;H_1)\cap L^{2}(0,T;V_1)$\\
	        \item $u,\theta$ satisfy, $\forall v\in V_0$ and $\forall\zeta\in V_1$,  a.e.t 
	    \begin{equation}
	    \begin{split}
	         &\frac{d}{dt}(u,v)_0+\nu((u,v))_0+(B_0(u,u),v)=(\theta\bm{e}_3,v)_0\\
	         &\frac{d}{dt}(\theta,\zeta)_1+\kappa ((\theta,\zeta))_1+(B_1(u,\theta),\zeta)_1- (u.\bm{e}_3,\zeta)_1=0.
	    \end{split}
	    \end{equation}
	    \end{itemize}
	    A Leray-Hopf weak solution additionally satisfies, a.e. $s$, and for all  $t \ge s$,  the energy inequality
	    \begin{align}
	         |u(t)|^2+\int_s^t\nu\|u(\sigma)\|^2d\sigma&\leq |u(s)|^2+\int_s^t(\theta \bm{e}_3(\sigma),u(\sigma))d\sigma.\label{lerener1}\\
	        |\theta(t)|^2+\int_s^t\kappa\|\theta(\sigma)\|^2d\sigma&\leq |\theta(s)|^2+\int_s^t(u\cdot \bm{e}_3(\sigma),\theta(\sigma))d\sigma\label{lerener2}.
	    \end{align}
	\end{definition}
	\subsection{Interpolant Operators}
	A finite rank, bounded linear operator $I_h : L^2(\Omega)\to L^2(\Omega)$ is said to be a type-\rom{1} interpolant observable if there exists a dimensionless constant $c>0$
	such that
	\begin{equation}\label{intest}
	|I_h(v)|\leq c|v|~~ \forall v\in L^2(\Omega) \text{ and } |I_h(v)-v| \leq ch\|v\|~~ \forall v\in H^1(\Omega).
	\end{equation}
	We look at two major examples of type-\rom{1} interpolants.
	\begin{itemize}
	    \item \textbf{Modal interpolation}: In this case $I_h u=P_K^0(u)$ with $h\sim 1/K$, where $P_K^0$ denotes the orthogonal projection onto the space spanned by the first $K$ eigenvectors of the Stokes operator$A_0.$  Indeed, one can easily check that it satisfies \eqref{intest}:
	\begin{equation}\label{modalest}
	|P_K(v)|\leq |v|~~\forall v\in L^2(\Omega) ~~\text{ and } ~~|P_K(v)-v| \leq\frac{1}{\lambda_K^{1/2}}\|v\|~~\forall v\in H^1(\Omega).
	\end{equation}
	where $\lambda_K=\min\{\lambda_K^0,\lambda_K^1\}.$\\
	\item \textbf{Volume interpolation}: In this case, $\Omega$ is partitioned into $N$ smaller cuboids $Q_{\alpha}$, where $\alpha\in\mathcal{J}=\left\{(j, k, l)\in \mathbb{N}\times \mathbb{N}\times \mathbb{N}:1\leq j, k, l\leq \sqrt[3]{N}\right\}$. Each cuboid is of diameter  $h=\sqrt{2L^2+1}/\sqrt[3]{N}$. The interpolation operator is defined as follows:
	\begin{equation}\label{volint}
	I_h(v) = \sum_{\alpha\in\mathcal{J}}\bar{v}_{\alpha}\chi_{Q_{\alpha}}(x),~~ v\in L^2(\Omega) 
	\end{equation}
	where 
	\begin{equation}\label{valpha}
	\bar{v}_{\alpha} = \frac{1}{|Q_\alpha|}\int_{Q_{\alpha}}v(x)dx,
	\end{equation}
	and $|Q_\alpha|$ is the volume of $Q_\alpha$. For $v\in V_0$, we define $I_h(v)=\left(I_h(v_1),I_h(v_2),I_h(v_3)\right),$ where $v=(v_1,v_2,v_3)$.
	\end{itemize}
 \subsection{Well-Posedness}
    	\begin{definition}\label{definesol}
	    $(u,\theta)$ is said to be a weak solution to \eqref{main1}-\eqref{main3} if for all $T>0$,
	    \begin{itemize}
	        \item $u\in L^{\infty}(0,T;H_0)\cap L^{2}(0,T;V_0)$\\
	        \item $\theta\in L^{\infty}(0,T;H_1)\cap L^{2}(0,T;V_1)$\\
	        \item $u,\theta$ satisfy, $\forall v\in V_0$ and $\forall\zeta\in V_1$,  a.e.t 
	    \begin{equation}
	    \begin{split}
	         &\frac{d}{dt}(u,v)_0+\nu((u,v))_0+(B_0(u,u),v)=(\theta\bm{e}_3,v)_0\\
	         &\frac{d}{dt}(\theta,\zeta)_1+\kappa ((\theta,\zeta))_1+(B_1(u,\theta),\zeta)_1- (u.\bm{e}_3,\zeta)_1=0.
	    \end{split}
	    \end{equation}
	    \end{itemize}
	    A weak solution is said to be a strong/regular solution if it also belongs to $$\left[L^{\infty}(0,T; V)\cap L^2(0,T;D(A))\right]\times\left[L^{\infty}(0,T; V_1)\cap L^2(0,T;D(A_1))\right].$$
	\end{definition}
	
	\begin{remark}
	    	Similar to the case of 3-D Navier-Stokes, we have existence but not uniqueness of weak solution. Also, given initial data in $V_0\times V_1$, a unique strong solution exists for some $[0,T]$ and a strong solution is unique in the larger class of Leray-Hopf weak solutions.
	\end{remark}

	\subsection{Existence Of Weak Solution}

	In this section, we prove the existence of a weak solution to the data assimilated Boussinesq equation.
	Our data assimilated algorithm is given by the solution $(w,\eta)$ of 
	\begin{align}
	&\frac{dw}{dt}+\nu A_0(w)+B_0(w,w)=P_{\sigma}(\eta \bm{e}_3)+\mu(P_{\sigma}(I_h(u)-I_h(w))) \label{main4}\\ 
	&\frac{d\eta}{dt}+\kappa A_1(\eta)+B_1(w,\eta)- w.\bm{e}_3=0\label{main5}\\
	&w(x,0)=0, \eta(x,0)=0.\label{six}
	\end{align} 
	The Galerkin approximation $(w_n,\eta_n)$ of $(w,\eta)$ satisfies
	\begin{align}
	&\frac{dw_n}{dt}+\nu A_0(w_n)+P_nB_0(w_n,w_n)=(\eta_n \bm{e}_3)+\mu(P_{n}(I_h(u)-I_h(w_n))) \label{main6}\\ 
	&\frac{d\eta_n}{dt}+\kappa A_1(\eta_n)+P_n B_1(w_n,\eta_n)- w_n\cdot\bm{e}_3=0\label{main7}\\
	&w_n(x,0)=0, \eta_n(x,0)=0 \label{main8}
	\end{align} 
	where, by abuse of notation, $P_n$ denotes the projection onto the space spanned by the first $n$ eigenvectors of $A_i$ for $i=0,1$.
	\begin{theorem}\label{weakexist}
			Let $(u,\theta)$ be a Leray-Hopf weak solution to \eqref{main1}-\eqref{main3} for  all $t\geq 0$ and $I_h$ be any type 1 interpolant. Let $h_0>0$ be such that
		\[
			h_0^2= \frac{\nu\kappa\lambda_1}{16c}.
			\]
		Then, provided $h\leq h_0$ and $\mu$ is chosen satisfying
		\begin{equation}\label{weakcondn}
		\frac{\nu}{2ch_0^2}\leq\mu\leq \frac{\nu}{2ch^2}\quad \left(\frac{1}{\lambda_1}\sim L^2\right),
		\end{equation}
		there exists a weak solution $(w,\eta)$ of \eqref{main4} such that for any $T>0$,
		$$w\in L^{\infty}(0,T;H_0)\cap L^{2}(0,T;V_0), \text{ and } \eta\in L^{\infty}(0,T;H_1)\cap L^{2}(0,T;V_1).$$
	\end{theorem}
	\begin{proof}
		Taking the inner product of \eqref{main6} with $w_n$ and \eqref{main7} with $\eta_n$ and adding, we obtain
		\begin{equation}
		\begin{split}
		\frac{1}{2}\frac{d}{dt}\left(|w_n|_0^2+|\eta_n|_1^2\right)+\nu\|w_n\|_0^2+\kappa\|\eta_n\|_1^2&\leq 2|\eta_n|_1|w_n|_0+\mu(I_h(u)-I_h(w_n),w_n)\\
		&\leq 2|\eta_n|_1|w_n|_0+\mu(w_n-I_h(w_n),w_n)\\
		&-\mu|w_n|_0^2+\mu(I_h(u),w_n)
		\end{split}
		\end{equation}
		We bound each of the terms on the RHS below.\\
		First, using Young’s inequality, \eqref{weakcondn} and \eqref{poincare}, we have
		\begin{equation}
		\begin{split}
		2|\eta_n|_1|w_n|_0&\leq \frac{2|w_n|_0^2}{\kappa\lambda_1}+\frac{\kappa\lambda_1|\eta_n|_1^2}{2}\\
		&\leq \frac{\mu|w_n|_0^2}{4}+\frac{\kappa\|\eta_n\|_1^2}{2}.
		\end{split}
		\end{equation}
		Next, using \eqref{intest}, Cauchy-Schwartz, Young’s inequality and the second inequality in \eqref{weakcondn}, we obtain
		\begin{equation}\label{b1}
		\begin{split}
		\mu|(w_n-I_h(w_n),w_n)|_0&\leq \mu|w_n-I_h(w_n)|_0|w_n|_0\\
		&\leq \mu ch^2\|w_n\|_0^2+\frac{\mu}{4}|w_n|_0^2\\
		&\leq \frac{\nu}{2}\|w_n\|_0^2+\frac{\mu}{4}|w_n|_0^2,
		\end{split}
		\end{equation} 
		Lastly, using Cauchy-Schwartz and Young’s inequality, we have
		\begin{equation}\label{b3}
		\begin{split}
		\mu|(I_h(u),w_n)|_0&\leq \mu|I_h(u)|_0^2+\frac{\mu}{4}|w_n|_0^2.
		\end{split}
		\end{equation}
		Combining all the estimates, we obtain
		\begin{align}\label{tbound}
		\frac{d}{dt}\left(|w_n|_0^2+|\eta_n|_1^2\right)+\nu\|w_n\|_0^2+\kappa\|\eta_n\|_1^2+\frac{\mu}{2}|w_n|_0^2&\leq 2\mu|I_h(u)|_0^2.
		\end{align}
		Splitting the above inequality and using \eqref{poincare}, we have
		\begin{equation*}
		\begin{split}
			\frac{d}{dt}|w_n|_0^2+\frac{\mu}{2}|w_n|_0^2\leq 2\mu|I_h(u)|_0^2
		\end{split}
		\end{equation*}
		and
		\begin{align*}
		\frac{d}{dt}|\eta_n|_1^2+\kappa\lambda_1|\eta_n|_1^2\leq 2\mu|I_h(u)|_0^2.
		\end{align*}
		Applying Gronwall, \eqref{unibound}, \eqref{intest} and \eqref{main8}, we obtain
		\begin{equation}\label{weakbound}
		|w_n|_0^2\leq 4M_0^2,\text{ and }|\eta_n|_1^2\leq \frac{2\mu M_0^2}{\kappa\lambda_1}.
		\end{equation}
		Dropping the first and the last term on LHS of \eqref{tbound} and integrating on the interval $[s,s+1]$, we obtain
		\begin{equation}
		    \int_s^{s+1}\left(\nu\|w_n\|_0^2+\kappa\|\eta_n\|_1^2\right)\leq 2 \mu M_0^2.
		\end{equation}
		The remainder of the proof is similar to the proof of existence of weak solutions of the 3D NSE.
	\end{proof}
\begin{remark}
    The data assimilated solution we have obtained is also unique. Later, from Theorem \ref{W++}, we obtain uniqueness of solution on the time interval $[0.\infty)$. The same argument can be used to show uniqueness here on the time interval $[0,T]$.
\end{remark}

\subsection{Time Independent Bound On Data Assimilated Temperature}
The following theorem will establish a time independent bound on $\|\eta\|_{L^{2p}}$, which will be a stepping stone in proving convergence of the data assimilated solution to the actual solution.
\begin{theorem}\label{etan}
	Let $p\in \mathbb{N}$ and $\eta$ and $\eta_n$ be as in \eqref{main5} and \eqref{main8} respectively. Assume that the hypotheses of Theorem~\ref{weakexist} hold. Then $\|\eta\|_{L^{2p}}$ is uniformly bounded in time.
\end{theorem}
\begin{proof}
	Taking the inner product of \eqref{main7} with $\eta_n^{2p-1}$, we obtain
	\begin{equation}
	\left(\frac{d\eta_n}{dt},\eta_n^{2p-1}\right)_1-\kappa\left(\Delta\eta_n,\eta_n^{2p-1}\right)_1+\left((w_n\cdot\nabla)\eta_n,\eta_n^{2p-1}\right)_1=(w_n\cdot\bm{e}_3,\eta_n^{2p-1})_1
	\end{equation}
	We estimate each term below.
	\begin{equation}
	\begin{split}
	\frac{1}{2p}\frac{d\eta_n^{2p}}{dt}&=\eta_n^{2p-1}\frac{d\eta_n}{dt}\\
	\Rightarrow \left(\frac{d\eta_n}{dt},\eta_n^{2p-1}\right)_1&=\frac{1}{2p}\frac{d}{dt}\int_{\Omega}\eta_n^{2p}dx\\
	&=\frac{1}{2p}\frac{d}{dt}\|\eta_n\|^{2p}_{L^{2p}}
	\end{split}
	\end{equation}
	Integrating by part, we obtain
	\begin{equation}
	\begin{split}
	-\kappa(\Delta\eta_n,\eta_n^{2p-1})_1&=(2p-1)\kappa\left(\nabla\eta_n,\eta_n^{2p-2}\nabla\eta_n\right)_1\\
	&=(2p-1)\kappa\left(\eta_n^{p-1}\nabla\eta,\eta_n^{p-1}\nabla\eta_n\right)_1\\
	&=\frac{2p-1}{p^2}|\nabla\left(\eta_n^p\right)|^2_1
	\end{split}
	\end{equation}
	\begin{equation}
	\begin{split}
	\left((w_n\cdot\nabla)\eta_n,\eta_n^{2p-1}\right)_1&=-\left(\left(w_n\cdot\nabla\right)\eta_n^{2p-1},\eta_n\right)_1\\
	&=-(2p-1)\left(\eta_n^{2p-2}\left(w_n\cdot\nabla\right)\eta_n,\eta_n\right)_1\\
	&=-\frac{2p-1}{p^2}\left(\left(w_n\cdot\nabla\eta_n\right),\eta_n^{2p-1}\right)_1
	\end{split}
	\end{equation}
	Therefore, $\left(\left(w\cdot\nabla\eta_n,\eta_n^{2p-1}\right)\right)_1=0$.\\
	From Holder's inequality, we have
	\begin{equation}
	\begin{split}
	|(w_n\cdot\bm{e}_3,\eta_n^{2p-1})_1|\leq \|w_n\|_{L^{q_1}}\cdot\left\|\left(\eta_n^p\right)^{2-\frac{1}{p}}\right\|_{L^{q_2}}
	\end{split}
	\end{equation}
	where $q_1=\frac{6p}{4p+1}$ and $q_2=\frac{6p}{2p-1}$. Note that $q_1\leq 2~\forall p\in\mathbb{N}$ and since $\Omega$ is a bounded domain of finite measure, we obtain
	\begin{equation}
	\begin{split}
	|(w_n\cdot\bm{e}_3,\eta_n^{2p-1})_1|&\leq C|w_n|_0\cdot\left\|\left(\eta_n^p\right)^{2-\frac{1}{p}}\right\|_{L^{q_2}}\\
	&\leq C|w_n|_0\cdot\left(\int_{\Omega}\left(\eta_n^p\right)^6\right)^{1/q_2}\\
	&\leq C|w_n|_0\cdot\big(\left\|\eta_n^p\right\|_{L^{6}}\big)^\frac{2p-1}{p}\\
	\end{split}
	\end{equation}
	Applying the Sobolev embedding $L^6(\Omega)\hookrightarrow H^1$,Young's inequality and \eqref{weakbound}, we obtain
	\begin{equation}
	\begin{split}
	|(w_n\cdot\bm{e}_3,\eta_n^{2p-1})_1|&\leq C|w_n|_0\cdot\big(\left|\nabla\left(\eta_n^p\right)\right|_{1}\big)^\frac{2p-1}{p}\\
	&\leq C p^{2p-2}|w_n|_0^{2p}+\left(\frac{2p-1}{2p^2}\right)\left|\nabla\left(\eta_n^p\right)\right|^2_1\\
	&\leq C p^{2p-2} M_0^{2p}+\left(\frac{2p-1}{2p^2}\right)\left|\nabla\left(\eta_n^p\right)\right|^2_1
	\end{split}
	\end{equation}
	Combining all the estimates we obtain 
	\begin{eqnarray}
	\frac{1}{2p}\frac{d}{dt}\|\eta_n\|^{2p}_{L^{2p}}+\left(\frac{2p-1}{2p^2}\right)\left|\nabla\left(\eta_n^p\right)\right|^2_1\leq C p^{2p-2} M_0^{2p}
	\end{eqnarray}
	Substituting $\zeta_n=\eta^p_n$, we obtain
	\begin{eqnarray}
	\frac{1}{2p}\frac{d}{dt}|\zeta_n|^{2}_{1}+\left(\frac{2p-1}{2p^2}\right)\left|\nabla\zeta_n\right|^2_1\leq C p^{2p-2} M_0^{2p}
	\end{eqnarray}
	Since $\zeta_n=\eta^p_n$, it shares the same boundary conditions as $\eta_n$, and hence \eqref{poincare} is applicable to the above equation, giving us
	 \begin{eqnarray}
	 \frac{d}{dt}|\zeta_n|^{2}_{1}+\left(\frac{(2p-1)\lambda_1}{p}\right)\left|\zeta_n\right|^2_1\leq C p^{2p-1} M_0^{2p}
	 \end{eqnarray}
	 Now, applying Gronwall, we obtain
	 \begin{equation}
	\|\eta_n\|_{L^{2p}}^{2p}=|\zeta_n|_1^2\leq C \left(\frac{p^{2p}}{(2p-1)\lambda_1}\right) M_0^{2p}
	 \end{equation}
	 Therefore, we have
	 \begin{equation}\label{tempboundn}
	 \|\eta_n\|_{L^{2p}}\leq  \frac{Cp M_0 }{\left((2p-1)\lambda_1\right)^{\frac{1}{2p}}}\coloneqq S_{p,u}=S_{p}
	 \end{equation}
	 Hence $\eta_n$ is a bounded sequence in $L^{\infty}(0,T;L^{2p}(\Omega))$ and there exists a subsequence $\eta_{n_k}$ that converges to $\eta^*$ in $L^{\infty}(0,T;L^{2p}(\Omega))$ in the weak star topology. In other words
	 \begin{equation}\label{yayayay}
	     \langle\eta_{n_k},\xi\rangle_{X',X}\to\langle\eta^*,\xi\rangle_{X',X}~~~\forall \xi\in X=L^{1}(0,T;L^q(\Omega)),
	 \end{equation}
	 where $2p$ and $q$ are Holder conjugates. We also know that $\eta_n$, and hence $\eta_{n,k}$, converges to $\eta$ in $L^{\infty}(0,T;L^{2}(\Omega))$ in the weak star topology. In other words
	 \begin{equation}
	     \langle\eta_{n_k},\xi\rangle_{X_1',X_1}\to\langle\eta,\xi\rangle_{X_1',X_1}~~~\forall \xi\in X_1=L^{1}(0,T;L^2(\Omega)).
	 \end{equation}
	 Now since $p\geq 1$ and $\Omega$ is bounded, $L^{2p}(\Omega)\subset L^{2}(\Omega)$ and $L^{2}(\Omega)\subset L^{q}(\Omega)$. This would mean that 
	 \begin{equation}
	     L^{\infty}(0,T;L^{2p}(\Omega))\text{ ~~is a subspace of ~~} L^{\infty}(0,T;L^{2}(\Omega))
	 \end{equation}
     and the convergence in \eqref{yayayay} will also hold in the weaker topology, giving us
     \begin{equation}
	     \langle\eta_{n_k},\xi\rangle_{X_1',X_1}\to\langle\eta^*,\xi\rangle_{X_1',X_1}~~~\forall \xi\in X_1=L^{1}(0,T;L^2(\Omega)).
	 \end{equation}
	 This would mean the sequence $\eta_{n_k}$ has two weak star limits in $L^{1}(0,T;L^2(\Omega))$. Hence $\eta^*=\eta$ and 
	 \begin{equation}\label{tempbound}
	 \|\eta\|_{L^{2p}}\leq  \frac{Cp M_0 }{\left((2p-1)\lambda_1\right)^{\frac{1}{2p}}}\coloneqq S_{p,u}=S_{p}
	 \end{equation}
\end{proof}

\subsection{Global Existence Of Strong Solution}
Another key result we will need in order to show convergence of the data assimilated solution to the actual solution is the regularity of the data assimilated velocity $w$. In order to do so, we will have to impose conditions on our data and look at the term $\|I_h(u)\|_0$.
For the case of modal interpolation, in addition to satisfying \eqref{intest}, $I_h$ also satisfies
\begin{equation}\label{mod}
\|I_h(v)\| \leq \|P_N(u)\|\leq c\|v\|~~ \forall v\in H^1(\Omega),~N\sim\frac{1}{h}
\end{equation}
The piece-wise constant volume interpolant as defined in \eqref{volint} may not satisfy \eqref{mod} due to the lack of regularity of the characteristic function. In order to establish an inequality similar to \eqref{mod} and also an inequality explicitly in terms of the data for volume interpolation, we define a smoothed volume interpolant operator $\tilde{I}_h$ that satisfies 
\begin{equation}\label{moddat}
\|\tilde{I}_h(v)\|_0^2 \leq Ch\sum_{\alpha\in\mathcal{J}}|\bar{v}_{\alpha}|^2 \leq C\|v\|_0^2~~ \forall v\in V_0,
\end{equation}
where $\bar{v}_\alpha$ is as in \eqref{valpha}. The proof of \eqref{moddat} and the justification as to why $\tilde{I}_h$ is a type-\rom{1} interpolant is provided in the appendix.\\
Hence, we can now define a modified type-\rom{1} interpolant which encompasses the usual modal interpolant and the smoothed volume interpolant.

    Next, for a general Leray-Hopf weak solution $u$ of \eqref{main1}-\eqref{main3}, we define the quantity $M_{h,u}$($=M_{h,u,T}$) as
     \begin{equation}\label{mh}
       \begin{split}
        M_{h,u}^2=32\sup_{0\leq t\leq T}\begin{dcases}
        \displaystyle\|P_N(u)\|^2\sim \sum_{|k|\leq N}|\lambda_{k,0}|^2|\hat{u}(k)|^2,~N\sim\frac{1}{h}&(\text{Modal})\\[10pt]
         \displaystyle Ch \sum_{\alpha}|\bar{u}_{\alpha}|^2,
          ~\bar{u}_\alpha = \frac{1}{|Q_\alpha|}\int_{Q_\alpha} u &(Volume)
        \end{dcases}
        \end{split}
        \end{equation}
    where, $\lambda_{k,0}$ is the $k^{th}$ smallest eigenvalue of $A_0$ corresponding to the eigenvector $\phi_{k,0}$ and $\displaystyle P_N(u)=\sum_{k=1}^N\hat{u}(k)\phi_{k,0}$.

\begin{remark}\label{finMh}
    Observe that $M_{h,u}$, as defined in \eqref{mh}, depends only on the data and is always finite. In the volume interpolation case, from \eqref{volint}, we have
    \begin{equation}
        \begin{split}
            \sum_{\alpha\in\mathcal{J}}|\bar{u}_{\alpha}|^2 \leq\sum_{\alpha\in\mathcal{J}}\left(\left|\frac{1}{|Q_\alpha|}\int_{Q_{\alpha}}u(x)dx\right|\right)^2&\leq \frac{N^2}{|\Omega|^2}\sum_{\alpha\in\mathcal{J}}\|u\|^2_{L^1(Q_{\alpha})}\\
            &\leq \frac{N}{|\Omega|}\sum_{\alpha\in\mathcal{J}}\|u\|^2_{L^2(Q_{\alpha})}\\
            &\leq \frac{N}{|\Omega|}|u|_0^2\\
            &\leq \frac{N}{|\Omega|}M_0^2\hspace{5mm}\text{(which is finite)}.
        \end{split}
    \end{equation}
    For the modal case, from \eqref{mh},we may write
    \begin{equation}
        \begin{split}
            \|P_N(u)\|^2\leq \lambda_N|u|_0\leq \lambda_N M_0~\text{(which is finite).}
        \end{split}
    \end{equation}
\end{remark}

\noindent
Next, we have the following theorem that establishes the regularity of $w$.
\begin{theorem}\label{velreg1}
	Let $\tilde{I}_h$ be a modified general type-\rom{1}  interpolant, $M_{h,u}$ be as in \eqref{mh} and $0<T\leq \infty$.
	Let $h_0>0$ be be given by
	\begin{equation}  \label{maincondition}
        h_0^{-2}=\frac{4c}{\nu}\max\left\{\frac{8}{\kappa\lambda_1},
	\frac{2}{\kappa}\left(1+\frac{2}{\lambda_1^2}\right)\right\}.
	\end{equation}
	Assume that for some  $0<h\leq h_0$
	\begin{equation}
	\frac{16cM_{h,u}^4}{\nu^3}\leq\frac{\nu}{4ch^2}.
	\end{equation}
	Let $\mu$ be chosen such that
	\begin{equation}\label{vstrong2}
	    \max\left\{\frac{\nu}{4ch_0^2},\frac{16cM_{h,u}^4}{\nu^3}\right\}\leq\mu\leq\frac{\nu}{4ch^2}.
	\end{equation}
	Then, the data assimilated fluid velocity is regular and satisfies $$\|w\|_0\leq M_{h,u}.$$
\end{theorem}
\begin{proof}
	Taking the inner product of \eqref{main6} with $A_0(w_n)$ and \eqref{main7} with $\eta_n$ and adding, we obtain
	\begin{align*}
	\frac{1}{2}\frac{d}{dt}\left(\|w_n\|_0^2+|\eta_n|_1^2\right) +\nu|A_0w_n|_0^2+\kappa\|\eta_n\|_1^2&=-(B_0(w_n,w_n),A_0w_n)_0+\mu(w_n-\tilde{I}_hw_n,A_0w_n)_0\\
	&-\mu\|w_n\|_0^2+\mu(\tilde{I}_hu,A_0w_n)_0\\
	&+(\eta_n,A_0(w_n)\cdot \bm{e}_3)_1+(\eta_n,w_n\cdot \bm{e}_3)_1
	\end{align*}
	First, applying \eqref{nolinest1} and Young’s inequality, we have
	\begin{align*}
	|(B_0(w_n,w_n),A_0w_n)|_0&\leq c\|w_n\|_0^{3/2}|A_0(w_n)|_0^{3/2}\\
	&\leq \frac{c}{\nu^3}\|w_n\|_0^6+\frac{\nu}{4}|A_0w_n|_0^2.
	\end{align*}
	Next, from \eqref{modalest}, Cauchy-Schwartz, Young’s inequality and the second inequality in \eqref{vstrong2}, we have
	\begin{align*}
	|\mu(w_n-\tilde{I}_hw_n,A_0w_n)|_0&\leq\mu ch\|w_n\|_0|A_0w_n|_0\\
	&\leq \frac{\mu^2ch^2}{\nu}\|w_n\|_0^2+\frac{\nu}{4}|A_0w_n|_0^2\\
	&\leq \frac{\mu}{4}\|w_n\|_0^2+\frac{\nu}{4}|A_0w_n|_0^2,
	\end{align*}
	Also, from \eqref{mod}, \eqref{moddat} and \eqref{mh}, we see that, for a modified type-\rom{1} interpolant $\tilde{I}_h$, we may write
    \begin{equation}\label{mhbound}
       32 \|\tilde{I}_hu(t)\|_0 \leq M_{h,u},~t\in[0,T].
    \end{equation}
	For the remaining terms, we apply Cauchy-Schwarz, Young’s inequality and \eqref{poincare} to obtain
	\begin{align*}
	\mu(\tilde{I}_h(u),A_0w_n)_0&\leq \mu\|\tilde{I}_hu(t)\|_0^2+\frac{\mu}{4}\|w_n\|_0^2,
	\end{align*}
	\begin{align*}
	\|\eta_n\|_1\|w_n\|_0&\leq \frac{\|\eta_n\|_1^2}{\mu}+\frac{\mu\|w_n\|_0^2}{4},
	\end{align*}
	\begin{align*}
	|\eta_n|_1|w_n|_0&\leq \frac{2|\eta_n|_1^2}{\mu\lambda_1}+\frac{\mu\lambda_1|w_n|_0^2}{8}\\
	&\leq\frac{2\|\eta_n\|_1^2}{\mu\lambda_1^2}+\frac{\mu\|w_n\|_0^2}{8}
	\end{align*}
	Combining estimates we obtain
	\begin{align*}
	\frac{1}{2}\frac{d}{dt}\left(\|w_n\|_0^2+|\eta_n|_1^2\right) +\frac{\nu}{2}|A_0w_n|_0^2+\left(\kappa-\frac{1}{\mu}-\frac{2}{\mu\lambda_1^2}\right)\|\eta_n\|_1^2&+\left(\frac{\mu}{8}-\frac{c}{\nu^3}\|w_n\|_0^4\right)\|w_n\|_0^2\\
	&\leq \mu\sup_{[0,T]}\|\tilde{I}_hu(t)\|_0^2.
	\end{align*}
	Let $[0, T_1]$ be the maximal interval on which $\|w_n(t)\|\leq M_{h,u}$ for $t\in[0, T_1]$, 
	where $M_{h,u}$ as in \eqref{moddat}. Note that $T_1 > 0$ exists because we have $w_n(0)=0$. Assume
	that $T_1<T$ . Then by continuity, we must have $\|w_n(T_1)\|=M_{h,u}$. Applying \eqref{unibound},\eqref{weakbound}, the first inequality in \eqref{vstrong2} and dropping all terms except the first and the last term on the LHS, we obtain
	\begin{align*}
	\frac{d}{dt}\|w_n\|_0^2 +\frac{\mu}{8}\|w_n\|_0^2\leq 2\mu\sup_{[0,T]}\|\tilde{I}_hu(t)\|_0^2.
	\end{align*}
		Applying Gronwall's inequality and \eqref{mhbound}, we obtain
	\begin{align*}
	\|w_n(t)\|_0^2 \leq 16\sup_{[0,T]}\|\tilde{I}_hu(t)\|_0^2\leq \frac{1}{2}M_{h,u}^{2}~~\forall t\in [0,T_1].
	\end{align*}
	This contradicts the fact that $\|w_n(T_1)\|=M_{h,u}$. Therefore $T_1\geq T$ and consequently, $\|w_n(t)\|\leq M_{h,u}$ for all $t\in[0, T ].$ Passing to the limit as $n\to\infty$, we obtain the desired conclusion for $w$.
\end{proof}

\subsection{Synchronization For General Type-\rom{1} Interpolant}

We will now show that the data assimilated solution approaches the actual solution(synchronization) for the general type-\rom{1} interpolant. To show synchronization, in addition to \eqref{lerener1} and \eqref{lerener2}, we will need an energy equality for the data assimilation equation. Taking the inner product of \eqref{main4} and \eqref{main5} with $w$ and $\eta$ respectively, we obtain
	\begin{align}
	    |w(t)|_0^2+2\nu\int_0^t\|w(s)\|_0^2ds&= |w(0)|_0^2 + 2\mu\int_0^t(I_h(u(s)-w(s)),w(s))_0ds\nonumber \\
	    &+2\int_0^t(w\cdot \bm{e}_3(s),\eta(s))_1ds.\label{datenergy1}\\
	     |\eta(t)|_1^2+2\kappa\int_0^t\|\eta(s)\|_1^2ds&= |\eta(0)|_1^2+2\int_0^t(w\cdot \bm{e}_3(s),\eta(s))_1ds.\label{datenergy2}
	\end{align}
	We now state the following lemma.
	\begin{lemma}\label{sather-serrin}
	    Let $(u,\theta)$ be the general Leray-Hopf weak solution of \eqref{main1}-\eqref{main3} and $(w,\eta)$ be the weak solution to \eqref{main4}-\eqref{six}. Let $\tilde{w}=w-u$, $\tilde{\eta}=\eta-\theta$ and $0<T\leq\infty.$ Then, 
	    \begin{align}
	         |\tilde{w}(t)|_0^2+2\nu\int_0^t\|\tilde{w}(s)\|_0^2ds &\leq |\tilde{w}(0)|_0^2+2\int_0^tB_0(\tilde{w}(s),\tilde{w}(s)),w(s))_0ds-2\mu\int_0^t\left(I_h(\tilde{w}(s)),\tilde{w}(s)\right)_0ds\nonumber \\
	    &+2\int_0^t(\tilde{w}(s)\cdot\bm{e}_3,\tilde{\eta}(s))_1ds.\label{ahah1}\\
	         |\tilde{\eta}(t)|_1^2+2\kappa\int_0^t\|\tilde{\eta}(s)\|_1^2ds &\leq |\tilde{\eta}(0)|_1^2+2\int_0^tB_1(\tilde{w}(s),\tilde{\eta}(s)),\eta(s))_1ds+2\int_0^t(\tilde{w}(s)\cdot\bm{e}_3,\tilde{\eta}(s))_1ds\label{ahah2}
	    \end{align}
	\end{lemma}
	\begin{proof}
	    Recall that we denote by $\phi^0_j$ the eigenfunction corresponding to $\lambda_j^0$, $j=1,2,\dots,$ with $\lambda^0_1\leq\lambda^0_2\dots$ being the eigenvalues of $A_0.$ Let $H_n^0$ denote the span of $\phi^0_1,\phi^0_2,\dots \phi^0_n$ and $P_n^0$ denote the projection onto the space $H_n^0$. Also, let $P_n^0(u)=u_n$ and $P_n^0(w)=w_n$. It is important to note that $u_n$ and $w_n$ in this lemma are projections of $u$ and $w$ onto a finite dimensional subspace and must not be confused with the Galerkin projections of $u$ and $w$. We observe that since $u_n$ finite dimensional, $\displaystyle\frac{\partial u_n}{dt}$ exists in the classical sense. We hence have
	    \begin{equation}\label{mixeps}
	    \begin{split}
	        \frac{d}{dt}(w_n,u_n)_0=\left(\frac{\partial w_n}{dt},u_n\right)_0+\left(w_n,\frac{\partial u_n}{dt}\right)_0
	        =&\left( -B_0(w,w)-\nu A_0(w)+\mu(I_h(u-w)+\eta\bm{e}_3),u_n\right)_0\\
	        &+\left( w_n,f-B_0(u,u)-\nu A_0(u)+\theta\bm{e}_3\right)_0.
	        \end{split}
	    \end{equation}
	    Using the fact that $P_n^0$ commutes with $A_0$, we obtain
	    \begin{equation}\label{commute}
	        (A_0(u),w_n)_0=(A_0(u),P_n(w))_0=(P_n(u),A_0(w)_0)=(u_n,A_0(w))_0.
	    \end{equation}
	    Integrating on the interval $[s,t]$, $0\leq s<t\leq T$ and applying \eqref{commute}, we obtain
	    \begin{equation}
	    \begin{split}
	         (w_n(t),u_n(t))_0&-(w_n(s),u_n(s))_0=\int_s^t\left((\eta(\sigma)\bm{e}_3,u_n(\sigma))_0+(\theta(\sigma)\bm{e}_3,w_n(\sigma))_0\right)d\sigma\\
	         &-2\nu\int_s^t(A_0(w)(\sigma),u_n(\sigma))_0d\sigma+2\mu\int_s^t(I_h(u(\sigma)-w(\sigma))),u_n(\sigma))_0d\sigma\\
	         &-\int_s^t\left\{\left(B_0(w(\sigma),w(\sigma)),u_n(\sigma)\right)_0+\left(B_0(u(\sigma),u(\sigma)),w_n(\sigma)\right)_0\right\}d\sigma.
	    \end{split}
	    \end{equation}
	    $u_n(\tau)\to u(\tau)$ in $H_0$ for each $\tau \in [0,T]$. Letting $n \to \infty$ and passing through the limit, we see that for almost all $s$ and $t$, $0<s<t<T:$
	    \begin{equation}\label{mixed}
	    \begin{split}
	         (w(t),u(t))&-(w(s),u(s))=\int_s^t\left((\eta(\sigma)\bm{e}_3,u(\sigma))_0+(\theta(\sigma)\bm{e}_3,w(\sigma))_0\right)d\sigma\\
	         &-2\nu\int_s^t(A_0(w)(\sigma),u(\sigma))_0d\sigma+2\mu\int_s^t(I_h(u(\sigma)-w(\sigma))),u(\sigma))_0d\sigma\\
	         &-\int_s^t\left\{\left(B_0(w(\sigma),w(\sigma)),u(\sigma)\right)_0+\left(B_0(u(\sigma),u(\sigma)),w(\sigma)\right)_0\right\}d\sigma
	    \end{split}
	    \end{equation}
	    Since $u$ is weakly continuous in $H_0$ and $w$ is strongly continuous in $H_0$, the function $$t\to(u(t),w(t))$$ is continuous and the relation \eqref{mixed} holds for all $s$ and $t$, $0\leq s<t\leq T.$ Moreover, we observe that 
	    $$\left(B_0(w,w),u\right)_0+\left(B_0(u,u),w\right)_0=\left(B_0(w-u,w),u\right)_0=\left(B_0(w-u,w-u),u\right)_0.$$
	    Incorporating this in \eqref{mixed} and letting $s=0$, we obtain
	    \begin{equation}\label{finalmix}
	    \begin{split}
	        (w(t),u(t))+2\nu\int_0^t((w,u))_0&=(w(0),u(0))_0+\int_s^t\left((\eta(\sigma)\bm{e}_3,u(\sigma))_0+(\theta(\sigma)\bm{e}_3,w(\sigma))_0\right)d\sigma\\
	        &-\int_0^t(B_0(\tilde{w}(s),\tilde{w}(s)),w(s))ds-2\mu\int_0^t I_h(\tilde{w}(s),u)ds
	    \end{split}
	    \end{equation}
	    Adding \eqref{lerener1} and \eqref{datenergy1} and subtracting two times \eqref{finalmix}, we obtain \eqref{ahah1}. Repeating the same arguments for $$\frac{d}{dt}(\eta_n,\theta_n)_1=\left(\frac{\partial \eta_n}{dt},\theta_n\right)_1+\left(\eta_n,\frac{\partial \theta_n}{dt}\right)_1,$$ we obtain \eqref{datenergy2}.
	\end{proof}
With the above lemma in place, we will now prove synchronization in the following theorem.
\begin{theorem}\label{sync}
	Let $I_h$ be the general type-\rom{1} interpolant, $0<T\leq \infty$ and $M_{h,u}$ as in \eqref{mh}.
	Also, let $h_0>0$ be defined as
	\begin{equation}\label{h02}
	    h_0^{-2}=\max\left\{\frac{32c}{\nu\kappa\lambda_1},
	\frac{8c}{\nu\kappa}\left(1+\frac{2}{\lambda_1^2}\right),\frac{64CS_2^8}{\nu^4\kappa^4}\right\},
	\end{equation}
	where $S_2$ is as in \eqref{tempbound} with $p=2$. We define $\tilde{w}=w-u$ and $\tilde{\eta}=\eta-\theta$. Assume that for some $h\leq h_0$,
	\begin{equation}
    \frac{16cM_{h,u}^4}{\nu^3}\leq\frac{\nu}{4ch^2}.
	\end{equation}
	Let $\mu$ be chosen such that
	\begin{equation}\label{vsync1}
        \max\left\{\frac{\nu}{4ch_0^2},\frac{16cM_{h,u}^4}{\nu^3}\right\}\leq\mu\leq\frac{\nu}{4ch^2}.
	\end{equation}
	Then,
	\begin{equation}
	\left(|\tilde{w}(t)|_0^2+|\tilde{\eta}(t)|_1^2\right)\leq \left(|\tilde{w}(0)|_0^2+|\tilde{\eta}(0)|_1^2\right)e^{-\alpha t}~~\forall t\in[0,T]
	\end{equation}
	where $\displaystyle\alpha=\min\left\{\frac{\mu}{4},\frac{\kappa\lambda_1}{2}\right\}$. In particular, if in the statement of Theorem~\ref{velreg1}, $T=\infty$, then
	\begin{equation}
	\lim_{t\to\infty}\left(|\tilde{w}(t)|_0^2+|\tilde{\eta}(t)|_1^2\right)=0.
	\end{equation}
\end{theorem}

\begin{proof}
	Adding equations \eqref{datenergy1} and \eqref{datenergy2}, we obtain
	\begin{equation}\label{vcon1}
	\begin{split}
	   |\tilde{w}(t)|_0^2 + |\tilde{\eta}(t)|_1^2+2\nu\int_0^t\|\tilde{w}(s)\|_0^2ds&+2\kappa\int_0^t\|\tilde{\eta}(s)\|_1^2ds\leq |\tilde{w}(0)|_0^2+|\tilde{\eta}(0)|_1^2-2\mu|\tilde{w}|^2\\
	   &+2\int_0^t(B_0(\tilde{w}(s),\tilde{w}(s)),w(s))_0ds+2\int_0^t(B_1(\tilde{w}(s),\tilde{\eta}(s)),\eta(s))_1ds\\
	   &+2\mu\int_0^t\left(\tilde{w}-I_h(\tilde{w}(s),\tilde{w}(s))\right)_0ds+4\int_0^t(\tilde{w}(s)\cdot\bm{e}_3,\tilde{\eta}(s))_1ds
	\end{split}
	\end{equation}
	We bound each term on the RHS.\\
	First, applying \eqref{nolinest2}, Cauchy-Schwartz and Young’s inequality, we obtain
	\begin{equation}
	\begin{split}
	|B_0(\tilde{w},\tilde{w}),w)_0|&\leq|(B_0(\tilde{w},w),\tilde{w}|)|_0\\
	&\leq c|\tilde{w}|_0^{1/2}\|\tilde{w}\|_0^{3/2}\|w\|_0\\
	&\leq \frac{c}{\nu^3}\|w\|_0^4|\tilde{w}|_0^2+\frac{\nu}{2}\|\tilde{w}\|_0^2.
	\end{split}
	\end{equation}
	Next, using Holder's inequality, \eqref{lady}, Young's inequality and \eqref{tempbound}, we have
	\begin{equation}
	\begin{split}
	|(B_1(\tilde{w},\tilde{\eta}),\eta)|_1&\leq \|\tilde{w}\|_{L^{4}}\|\nabla\tilde{\eta}_n\|_{L^2}\|\eta\|_{L^4}\\
	&\leq |\tilde{w}|_0^{1/4}\|\tilde{w}\|_0^{3/4}\|\tilde{\eta}\|_{1}\|\eta\|_{L^4}\\
	&\leq \frac{\nu}{4}\|\tilde{w}\|_0^2+\frac{\kappa}{2}\|\tilde{\eta}\|_{1}^2+\frac{C}{\nu^3\kappa^4}|\tilde{w}|_0^2\|\eta\|_{L^4}^8\\
	&\leq \frac{\nu}{4}\|\tilde{w}\|_0^2+\frac{\kappa}{2}\|\tilde{\eta}\|_{1}^2+\frac{C}{\nu^3\kappa^4}S_2^8|\tilde{w}|_0^2
	\end{split}
	\end{equation}
	Applying \eqref{intest}, Cauchy-Schwartz, Young’s inequality and \eqref{vstrong2}, we obtain
	\begin{equation}
	\begin{split}
	|\mu(\tilde{w}-I_h(\tilde{w}))|_0|\tilde{w}|_0&\leq \mu ch\|\tilde{w}\|_0|\tilde{w}|_0\\
	&\leq \mu ch^2\|\tilde{w}\|_0^2+\frac{\mu}{4}|\tilde{w}|_0^2\\
	&\leq \frac{\nu}{4}\|\tilde{w}\|_0^2+\frac{\mu}{4}|\tilde{w}|_0^2.
	\end{split}
	\end{equation}
	Using Cauchy-Schwartz, Young’s inequality and \eqref{poincare} , we obtain
	\begin{equation}
	\begin{split}
	4|\tilde{\eta}|_1|\tilde{w}|_0&\leq \frac{16|\tilde{\eta}|_1^2}{\mu}+\frac{\mu|\tilde{w}|_0^2}{4}\\
	&\leq \frac{16\|\tilde{\eta}\|_1^2}{\lambda_1\mu}+\frac{\mu|\tilde{w}|_0^2}{4}
	\end{split}
	\end{equation}

	Combining the above estimates and applying \eqref{poincare}, \eqref{h02},\eqref{vsync1} and the fact that $\|w\|_0\leq M_{h,u}$, we obtain

	\begin{equation}\label{qw}
	\begin{split}
	|\tilde{w}(t)|_0^2+|\tilde{\eta}(t)|_1^2+\frac{\mu}{4}\int_0^t|\tilde{w}(s)|_0^2ds+\frac{\kappa\lambda_1}{2}\int_0^t|\tilde{\eta}(s)|_1^2ds\leq |\tilde{w}(0)|_0^2+|\tilde{\eta}(0)|_1^2.
	\end{split}
	\end{equation}
	From \eqref{qw} we observe that
	\begin{equation}\label{qw1}
	     |\tilde{w}(t)|_0^2+|\tilde{\eta}(t)|_1^2\leq |\tilde{w}(0)|_0^2+|\tilde{\eta}(0)|_1^2.
	\end{equation}
	Let $$\displaystyle\alpha=\min\left\{\frac{\mu}{4},\frac{\kappa\lambda_1}{2}\right\}.$$ Changing the interval of integration from $[0,t]$ to $[t-1/\alpha,t]$, applying Lemma \ref{titi} with $y=|\tilde{w}(t)|^2+|\tilde{\eta}(t)|^2$ and using \eqref{qw1}, we obtain
\begin{equation}
	\left(|\tilde{w}(t)|^2+|\tilde{\eta}(t)|^2\right)\leq \left( |\tilde{w}(0)|_0^2+|\tilde{\eta}(0)|_1^2\right)e^{-\alpha\left(t-\frac{1}{\alpha}\right)}~~\forall t\in[0,T].
	\end{equation}

\end{proof}
\section{Forward Determining map}
We now describe the lifting map introduced in \cite{DM}, which played a crucial role in obtaining the determining form for evolution equations. In some sense, this map is proposed as a substitute, or alternative, for the notion of inertial manifold. Using the ideas in \cite{IHP}, we construct a similar map for the 3-D Boussinesqu(and later for the 3-D Navier-Stokes) equation. Notably, the conditions we require for the construction of this map here(like in \cite{IHP}) are weaker than the ones considered in \cite{DM2}. In particular, we do not require any condition on the time derivative of the input function $v =v(t)$ in \eqref{mainn4}, below. This may be useful within the context of data assimilation, where $v(t)$ represents the observed spatial coarse-mesh data, which is usually noisy.\\

We begin by introducing spaces that contain the domain and ranges of the map $W^+$ and by introducing the evolution equation which yields the definition of $W^+$.\\
We denote by $L^2_b(\mathbb{R}^+;D(A_0))$ the functions in $L^2(\mathbb{R}^+;D(A_0))$ which are translation bounded, i.e.,
\begin{equation}\label{transbound}
   \sup_{s\geq 0} \int_s^{s+1}|A_0(u(r))|_{0}^2dr<\infty.
\end{equation}
Similarly, $L^2_b(\mathbb{R}^+;V)$ denotes the functions in $L^2(\mathbb{R}^+;V)$ which satisfy $$\displaystyle \sup_{s\geq 0}\int_s^{s+1}\|u(r)\|_{0}^2dr<\infty.$$
Let $$Y_+=\left(C_b(\mathbb{R}_+;V_0)\cap L^2_b(\mathbb{R}_+;D(A_0))\right)\times C_b(\mathbb{R}_+;H_1)$$ and  $$Z_+=\left(C_b(\mathbb{R}_+;H_0)\cap L^2_b(\mathbb{R}_+;V_0)\right)\times C_b(\mathbb{R}_+;H_1)$$
where, $C_b(I;B)$ denotes the space of all bounded and continuous functions over the interval $I$ with values in the Banach space $B$.
$Y_+$ and $Z_+$ are Banach spaces with norms
\begin{equation}\label{ynorm}
    \|u\|_{Y_+}=\left\{\|u\|_0^2+\sup_{s\geq 0}\int_s^{s+1}|A_0(u(r))|_{0}^2dr\right\}^{\frac{1}{2}}.
\end{equation}
and
\begin{equation}
    \|u\|_{Z_+}=\left\{|u|_0^2+\sup_{s\geq 0}\int_s^{s+1}\|u(r)\|_{0}^2dr\right\}^{\frac{1}{2}}
\end{equation}
respectively.
Moreover, let $X$ be the Banach space $$X_+=C_b(\mathbb{R}^+;(H^1(\Omega))^3)$$ equipped with the norm 
\begin{equation}\label{xnorm}
    \|v\|_{X_+}=\sup_{s\geq 0}\|v(s)\|_0.
\end{equation}
The observed spatial coarse-mesh data is denoted by $v(t)$. For the purpose of data assimilation, we consider the case $v\in X_+$ with $\|v\|_{X_+}\leq\rho\sim h^{-1/2}$ for some $\rho>0$. We use $B_{X_+}(\rho)$ to denote the closed ball in $X_+$ of radius $\rho$ centered at $0$.
 We also define the Banach space $$P_+=C_b(\mathbb{R}_+;H_0)\times C_b(\mathbb{R}_+;H_1)$$ with the norm
\begin{equation}\label{p+norm}
    \|(u,\theta)\|_{P_+}=\left\{|u|_0^2+|\theta|_1^2\right\}^{\frac{1}{2}}
\end{equation}
For $\tau\in I\subset\mathbb{R^+}$ and a Banach space $B$, we define the time translation $\tau_{\sigma}:C(I;B)\to C(I;B)$ as
\begin{equation}
    \tau_{\sigma}(u(t))=u(t+\sigma).
\end{equation}
Now given $v\in B_{X_+}(\rho)$, we consider the following initial-value problem:
	\begin{align}
	&\frac{dw}{dt}+\nu A_0(w)+B_0(w,w)=P_{\sigma}(\eta \bm{e}_3)+\mu(P_{\sigma}(v-I_h(w))) \label{mainn4}\\ 
	&\frac{d\eta}{dt}+\kappa A_1(\eta)+B_1(w,\eta)- w.\bm{e}_3=0\\
	&w(x,0)=0, \eta(x,0)=0\label{mainn5}.
	\end{align} 
	The corresponding Galerkin approximation is given by
	\begin{align}
	&\frac{dw_n}{dt}+\nu A_0(w)+B_0(w_n,w_n)=P_{n}(\eta_n \bm{e}_3)+\mu(P_{n}(v-I_h(w))) \label{gmainn4}\\ 
	&\frac{d\eta_n}{dt}+\kappa A_1(\eta_n)+B_1(w,\eta_n)- w_n.\bm{e}_3=0\\
	&w_n(x,0)=0, \eta_n(x,0)=0\label{gmainn5}.
	\end{align}
	We can now define the map $W_+$.
	\begin{definition}
	    Let $\rho>0$. We define the forward determining map $W_+:B_{X_+}(\rho)\to Y_+$ as
	\begin{equation}\label{W+}
	    W_+(v)=(w,\eta).
	\end{equation}
	\end{definition}

	\begin{theorem}\label{W++}
	Let $I_h$ be a general type-\rom{1}  interpolant, $v\in B_{X_+}(\rho)$ for some $\rho>0\left(\displaystyle\rho\sim 1/\sqrt{h}\right)$, $0<T\leq \infty$ and 
	\begin{equation}\label{vstrongg1}
	M_{h,u}^2=32\sup_{[0,T]}\|v\|_0^2.
	\end{equation}
		Also, let $h_0>0$ be defined as
	\begin{equation}\label{h04}
	    h_0^{-2}=\max\left\{\frac{32c}{\nu\kappa\lambda_1},
	\frac{8c}{\nu\kappa}\left(1+\frac{2}{\lambda_1^2}\right),\frac{64CS_2^8}{\nu^4\kappa^4}\right\},
	\end{equation}
	where $c$, $h$ are as in \eqref{intest} and $S_2$ is as in \eqref{tempbound} with $p=2$.	Assume that for some $h\leq h_0$,
	\begin{equation}
    \frac{16cM_{h,u}^4}{\nu^3}\leq\frac{\nu}{4ch^2}.
	\end{equation}
     Let $\mu$ be chosen such that
 \begin{equation}\label{vsyncc1}
    \max\left\{\frac{\nu}{4ch_0^2},\frac{16cM_{h,u}^4}{\nu^3}\right\}\leq\mu\leq\frac{\nu}{4ch^2}.
	\end{equation}
 Then, there exists a weak solution $(w,\eta)$ to \eqref{mainn4}-\eqref{mainn5}. Moreover, the following statements hold true
	\begin{enumerate}[(a)]
	    \item $\|w\|_0\leq M_{h,u}$\label{w1}\\
	    \item $W_+:B_{X_+}(\rho)\to Z_+$ is Lipschitz continuous.\label{w2}\\
	    \item Let $(u,\theta)$ be a weak solution \eqref{main1}-\eqref{main3}. Assume that $|I_h(u)(s)-v(s)|_0\to 0$ as $s\to\infty$. Then $\|W_+(v)(s)-(u(s),\theta(s))\|_{P_+}\to 0$ as $s\to\infty$\label{w3}\\
	    \item Let $v_1,v_2\in B_{X_+}(\rho)$, Then, $W_+(v_1)=W_+(v_2)$ iff $P_{\sigma}(v_1-v_2)=P_{\sigma}(\bar{v})=0,$ where $\bar{v}=v_1-v_2$.\label{w4}\\
	    \item For every $\sigma\in \mathbb{R}^+$,
	    \begin{equation}
	        W_+\circ\tau_\sigma(v)=\tau_\sigma\circ W_+(v).
	    \end{equation} \label{w5}
	\end{enumerate}
\end{theorem}
	\begin{proof}
	    The existence of a weak solution $(w,\eta)$ and the proof of \ref{w1} is obtained by repeating the proof of Theorem \ref{weakexist}
        and Theorem \ref{velreg1}, after replacing $I_h(u)$ with $v.$\\
        
        To prove \ref{w2}, we define $\bar{w}_n=w_n^1-w_n^2$, $\bar{\eta}_n=\eta_n^1-\eta_n^2$ and $\bar{v}_n=P_n(\bar{v}_1-\bar{v}_2)$, where $(w_n^1,\eta_n^1)$ and $(w_n^2,\eta_n^2)$ are Galerkin approximations of $(w^1,\eta^1)$ and $(w^2,\eta^2)$ respectively and $(w^1,\eta^1)=W_+(v_1)$ and $(w^2,\eta^2)=W_+(v_2)$.\\
        $\bar{w}_n$ and $\bar{\eta}_n$ satisfy the equations
	\begin{equation}\label{vconn1}
	\begin{split}
	\frac{d\bar{w}_n}{dt}+\nu A_0(\bar{w}_n)+B_0(\bar{w}_n,w_n)+B_0(u_n,\bar{w}_n)-\bar{\eta}_n \bm{e}_3&=\mu P_n(\bar{v}_n-I_h(\bar{w}_n))\\
	&=\mu \bar{v}_n+\mu P_n(\bar{w}_n-I_h(\tilde{w}_n))\\
	&-\mu\bar{w}_n
	\end{split}
	\end{equation}
	\begin{equation}\label{vconn2}
	\begin{split}
	\frac{d\bar{\eta}_n}{dt}+\kappa A_1(\bar{\eta}_n)+B_1(\bar{w}_n,\eta_n)+B_1(u,\bar{\eta}_n)-\bar{w}_n\cdot \bm{e}_3&=0
	\end{split}
	\end{equation}
	Taking the inner product of $\bar{w}_n$ and $\bar{\eta}_n$ with \eqref{vconn1} and \eqref{vconn2} respectively, we obtain
	\begin{equation}\label{bar}
	\begin{split}
	\frac{1}{2}\frac{d}{dt}\left(|\bar{w}_n|_0^2+|\bar{\eta}_n|_1^2\right)+\nu\|\bar{w}_n\|_0^2+\kappa\|\bar{\eta}_n\|_1^2+\mu|\bar{w}_n|_0^2&\leq |(B_0(\bar{w}_n,w_n),\bar{w}_n|)|_0+|(B_1(\bar{w}_n,\eta_n),\bar{\eta}_n)|_1\\
	&+2|\bar{w}_n|_0|\bar{\eta}_n|_1+|\mu(\bar{w}_n-I_h(\bar{w}_n))|_0|\bar{w}_n|_0\\
	&+\mu|\bar{v}_n|_0|\bar{w}_n|_0.
	\end{split}
	\end{equation}
	Applying Young's inequality, we obtain
	\begin{equation}
	\begin{split}
	\mu|\bar{v}_n|_0|\bar{w}_n|_0\leq \mu|\bar{v}_n|_0^2+\frac{\mu}{4}|\bar{w}_n|_0^2,
	\end{split}
	\end{equation}
	Bounding the other terms on the RHS of \eqref{bar} exactly as in Theorem \ref{sync}, we obtain
	\begin{equation}\label{combo}
	\begin{split}
	\frac{1}{2}\frac{d}{dt}\left(|\bar{w}_n|^2_0+\bar{\eta}_n|_1^2\right)+\nu\|\bar{w}_n\|_0+\left(\frac{\mu}{4}-\frac{c}{\nu^3}\|w_n\|_0^4-\frac{CS_2^8}{\nu^3\kappa^4}\right)& |\bar{w}_n|_0^2+\left(\frac{\kappa\lambda_1}{2}-\frac{4}{\mu}\right)|\bar{\eta}|_1^2\leq\mu|\bar{v}_n|_0^2.
	\end{split}
	\end{equation}
	We already have $\|w\|_0\leq M_{h,u}$. Now applying \eqref{vsyncc1} and \eqref{h04}, we have
	\begin{equation}
	\begin{split}
	\frac{d}{dt}\left(|\bar{w}_n|^2|_0+|\bar{\eta}_n|_1^2\right)+\frac{\mu}{4}|\bar{w}_n|_0^2+\frac{\kappa\lambda_1}{2}|\bar{\eta}_n|_1^2\leq 2\mu|\bar{v}_n|_0^2.
	\end{split}
	\end{equation}
	Next, applying Gronwall (for $t,\sigma\in[0,T]$ with $t>\sigma$) separately to $|\bar{w}_n|^2$ and $|\bar{\eta}_n|^2$, we obtain
	\begin{equation}
	\begin{split}
	|\bar{w}_n(t)|_0^2&\leq |\bar{w}_n(\sigma)|_0^2e^{-(\mu/4)(t-\sigma)}+8\sup_{t\in[0,T]}|\bar{v}_n|_0^2,\\
	|\bar{\eta}_n(t)|_1^2&\leq |\bar{\eta}_n(\sigma)|_1^2e^{-\left(\frac{\kappa\lambda_1}{2}\right)(t-\sigma)}+\frac{4\mu}{\kappa\lambda_1}\sup_{t\in\mathbb{R}_+}|\bar{v}_n|_0^2.
	\end{split}
	\end{equation}
	$\bar{w}_n$ and $\bar{\eta}_n $ converge respectively to $\bar{w}$ and $\bar{\eta}$ weakly. Also $\bar{v_n}\to\bar{v_n}$ in $B_{x_+}(\rho)$. Therefore, letting $n\to\infty$, setting $\sigma=0$ and noting that $\bar{w}(0)=\bar{\eta}(0)=0$, we obtain
	\begin{equation}\label{n1}
	\begin{split}
	|\bar{w}(t)|^2|_0&\leq 8\sup_{t\in\mathbb{R}_+}|\bar{v}|_0^2\\
	|\bar{\eta}(t)|_1^2&\leq \frac{4\mu}{\kappa\lambda_1}\sup_{t\in\mathbb{R}_+}|\bar{v}|_0^2.
	\end{split}
	\end{equation}
	Using \eqref{combo}, we may write
	\begin{equation}
	  \nu\|\bar{w}_n\|_0\leq\mu\sup_{t\in\mathbb{R}_+}|\bar{v}_n(t)|_0^2.
	\end{equation}
	Letting $n\to\infty$ and integrating both sides on the interval $[t,t+1]$, we obtain
	\begin{equation}\label{n2}
	 \sup_{t\in\mathbb{R}_+} \int_{t}^{t+1}\nu\|\bar{w}(s)\|_0ds\leq\mu\sup_{t\in\mathbb{R}_+}|\bar{v}(t)|_0^2.
	\end{equation}
	\eqref{n1} and \eqref{n2} together prove \ref{w2}.\\
	
	To prove \ref{w3}, we define $\tilde{w}=w-u$ and $\tilde{\eta}=\eta-\theta$, where $(w,\eta)=W_+(v).$
	Replacing $I_h(u)$ with $v$ in \eqref{datenergy1} and repeating the arguments of Lemma \eqref{sather-serrin}, we obtain a modified version of \eqref{ahah1}, given by
	\begin{equation}
	\begin{split}
	    |\tilde{w}(t)|_0^2+2\nu\int_{\sigma}^t\|\tilde{w}(s)\|_0^2ds \leq &|\tilde{w}(0)|_0^2+2\int_{\sigma}^tB_0(\tilde{w}(s),\tilde{w}(s)),w(s))_0ds\\
	    &+2\mu\int_{\sigma}^t\left(v(s)-I_hw(s),\tilde{w}(s)\right)_0ds+2\int_{\sigma}^t(\tilde{w}(s)\cdot\bm{e}_3,\tilde{\eta}(s))_1ds.\\
	    \leq &|\tilde{w}(0)|_0^2+2\int_{\sigma}^tB_0(\tilde{w}(s),\tilde{w}(s)),w(s))_0ds\\
	    &+2\mu\int_{\sigma}^t\left(v(s)-I_hu(s),\tilde{w}(s)\right)_0ds
	    +2\mu\int_{\sigma}^t\left(\tilde{w}-I_h\tilde{w}(s),\tilde{w}(s)\right)_0ds\\
	    &-\mu\int_{\sigma}^t\|\tilde{w}\|_0ds+2\int_{\sigma}^t(\tilde{w}(s)\cdot\bm{e}_3,\tilde{\eta}(s))_1ds..
	\end{split}
	\end{equation}
	Applying Cauchy-Schwarz and Young's inequality, we obtain
	\begin{equation}
	    2\mu\left|\left(v(s)-I_hu(s),\tilde{w}(s)\right)_0\right|\leq 2\mu|v(s)-I_hu(s)|_0^2+\frac{\mu}{2}|\tilde{w}(s)|_0^2.
	\end{equation}
	Also,
	\begin{equation}
	    |v|_0^2\leq\frac{1}{\lambda_1}\|v\|_0^2\leq \frac{1}{\lambda_1}\rho^2.
	\end{equation}
	and $|I_hu|_0\leq|u|_0\leq M_0.$ Hence it makes sense to talk about the term
	\begin{equation*}
	    \sup_{0\leq s\leq \infty}|I_h(u)(s)-v(s)|_0.
	\end{equation*}
	Proceeding exactly as in Theorem \ref{sync}, we obtain
	\begin{equation}\label{ccombo}
	\begin{split}
	|\tilde{w}(t)|^2_0+|\tilde{\eta}(t)|_1^2+\frac{\mu}{4}\int_{\sigma}^t|\tilde{w}(s)|_0^2ds+\frac{\kappa\lambda_1}{2}|\tilde{\eta}|_1^2\leq |\tilde{w}(\sigma)|_0^2+|\tilde{\eta}(\sigma)|_1^2+2\mu\sup_{\sigma\leq s\leq t}|I_h(u)(s)-v(s)|_0.
	\end{split}
	\end{equation}
	Next, letting $$\displaystyle\alpha=\min\left\{\frac{\mu}{4},\frac{\kappa\lambda_1}{2}\right\}$$
	and applying Corollary \ref{1lemmaaa1}, \eqref{weakbound} and \eqref{unibound}, we obtain
\begin{equation}
	\left(|\tilde{w}(t)|^2_0+|\tilde{\eta}(t)|_1^2\right)\leq 16Me^{-\alpha (t/2)}+\frac{2\mu}{\alpha}\sup_{t/2\leq s\leq t}|I_h(u)(s)-v(s)|_0.
	\end{equation}
	Taking the lim sup as $t\to\infty$ and using the hypothesis that $|I_h(u)(s)-v(s)|_0\to 0$ as $s\to\infty$, we see that
	\begin{equation}
	\lim_{t\to\infty}\left(|\tilde{w}(t)|^2_0+|\tilde{\eta}(t)|_1^2\right)=0,
	\end{equation}
	proving \ref{w3}.\\
	To prove \ref{w4}, we note that $\bar{w}=w_1-w_2=0$ and $\bar{\eta}=\eta_1-\eta_2=0$, where $(w_1,\eta_1)=W_+(v_1)$ and $(w_2,\eta_2)=W_+(v_2)$. Since $w_1$ and $w_2$ are regular, the term $\bar{w}$ is differentiable a.e on $\mathbb{R}^+$. From \eqref{mainn4}, we may write
	\begin{equation}
	    	\frac{d\bar{w}}{dt}+\nu A_0(\bar{w})+B_0(\bar{w},w_1)+B_0(w_2,\bar{w})-\bar{\eta} \bm{e}_3=\mu P_\sigma(\bar{v}-I_h(\bar{w}))
	\end{equation}
	Letting $\bar{w}=0$ and $\bar{\eta}=0$, we obtain $P_\sigma(\bar{v})=0.$ If $P_\sigma(\bar{v})=0$, We obtain $W_+(v_1)=W_+(v_2)$ from the Lipschitz continuity of $W_+.$\\
	To show \ref{w5}, we observe that $\tau_\sigma\circ W_+(v)$ is a solution of \eqref{mainn4}-\eqref{mainn5} corresponding to $\tau_\sigma(v)$. From the Lipschitz property of the map $W_+$, we have uniqueness of solution. Hence $$ W_+\circ\tau_\sigma(v)=\tau_\sigma\circ W_+(v).$$
\end{proof}
	\begin{corollary}\label{needcor}
        Let $I_h$ be a general type-\rom{1}  interpolant and $u$ be a weak solution \eqref{main1}-\eqref{main3}. Let the hypothesis of Theorem \ref{W++} hold. Then $\|W_+(I_h(u))(s)-(u(s),\theta(s))\|_{P_+}\to 0$ as $s\to\infty$
	\end{corollary}
	\begin{proof}
	    Applying part \ref{w3} of Theorem \ref{W++} with $v=I_h(u)$, we obtain the statement of the corollary.
	\end{proof}
	\begin{remark}
	    The fact that $W^+$ is Lipschitz continuous means that the way we reconstruct our solution from the data is "stable". Lipschitz continuity in turn implies uniqueness of solutions as well. Also, Corollary \ref{needcor} says that the solution obtained from type-\rom{1} interpolation data(for appropriate $h$) asymptotically approaches the actual solution.
	\end{remark}
\begin{theorem}
     Let $(u_1,\theta_1)$ and $(u_2,\theta_2)$ be two restricted Leray-Hopf weak solutions and $M_{h,u_i}$(for $i=1,2)$ be as in \eqref{mh}. Also, let $h_0>0$ be defined as
	\begin{equation}\label{h05}
	    h_0^{-2}=\max\left\{\frac{32c}{\nu\kappa\lambda_1},
	\frac{8c}{\nu\kappa}\left(1+\frac{2}{\lambda_1^2}\right),\frac{64CS_2^8}{\nu^4\kappa^4}\right\},
	\end{equation}
	where $c$, $h$ are as in \eqref{intest} and $S_2$ is as in \eqref{tempbound} with $p=2$.	Assume that on $[0,\infty)$, for $i=1,2$, for some $h\leq h_0$,
	\begin{equation}\label{cccc}
    \frac{16cM_{h,u_i}^4}{\nu^3}\leq\frac{\nu}{4ch^2}.
	\end{equation}
	Let $\mu$ be chosen such that
		\begin{equation}\label{vsynccc1}
    \max\left\{\frac{\nu}{4ch_0^2},\frac{16cM_{h,u_i}^4}{\nu^3}\right\}\leq\mu\leq\frac{\nu}{4ch^2}.
	\end{equation}
	Then, if $$\displaystyle \lim_{t\to\infty}|I_h(u_1)(t)-I_h(u_2)(t)|_0=0$$ then $$\lim_{t\to\infty}|u_1(t)-u_2(t)|=0\text{ and }\lim_{t\to\infty}|\theta_1(t)-\theta_2(t)|=0.$$
\end{theorem}
\begin{proof}
Repeating the proof of part \ref{w3} with $(u,\theta)=(u_1,\theta_1)$, $(w,\eta)=(u_2,\theta_2)$ and $v=I_hu_2$, we obtain the statement of the theorem.
\end{proof}
\begin{remark}
    The above theorem shows that $h\leq h_0$ satisfying \eqref{cccc} is \textit{asymptotically determining}.
\end{remark}
\section{An observable regularity criterion  on the Weak attractor of the three-dimensional Navier-Stokes equations}\label{attractor}
In this section we will be looking at solutions to the data assimilated 3-D Navier-Stokes equation on the weak attractor. We discuss well-posedness of the data assimilated equation for time ranging over all real numbers as well as the question of uniqueness and regularity for 3-D Navier-Stokes equation on the weak attractor when the \say{low modes} are known.\\

\noindent
The 3-D incompressible Navier-Stokes equations (3D NSE) on a domain with time independent forcing (assumed for simplicity) is given by
\begin{equation}\label{3dnav}
\begin{split}
    \frac{\partial u}{\partial t}+(u\cdot\nabla)u-\Delta u+\nabla p&=f\\
    \nabla\cdot u&=0\\
\end{split}
\end{equation}
Applying $P_{\sigma}$ to \eqref{3dnav}, we obtain
\begin{equation}\label{3dnav1}
\begin{split}
    \frac{\partial u}{\partial t}+B_0(u,u)+A_0(u)&=f.\\
    \nabla\cdot u&=0\\
    \end{split}
\end{equation}
where, by abuse of notation, we denote $P_{\sigma}(f)$ by $f$.
Applying $P_n^0$ to \eqref{3dnav1}, we obtain the Galerkin approximation to the system, given by 

\begin{equation}\label{3dnav2}
\begin{split}
    \frac{\partial u_n}{\partial t}+B_0(u_n,u_n)+A_0(u_n)&=P_n^0f\\
     \nabla\cdot u_n&=0\\
    \end{split}
\end{equation}
\subsection{Well-Posedness}
We begin by providing definitions of weak and strong solutions.
\begin{definition}\label{weakdef}
A (Leray-Hopf) weak solution on a time interval $I=[0,T] \subset \mathbb{R}$ is defined as a function $u = u(t)$ on $I$ with values in $H_0$ and satisfying the following properties:
\begin{itemize}
    \item $u\in L^{\infty}(0,T;H_0)\cap L^{2}(0,T;V_0)\cap C(0,T:V_0')$\\
    \item $\displaystyle\frac{du}{dt}\in L^{4/3}(0,T;V_0')$\\
    \item $u$ satisfies the functional equation \eqref{3dnav1} in the distribution sense on $I$, with values in $V_0'$;\\
    \item For almost all $t'\in I$, $u$ satisfies the following energy inequality
    \begin{equation}\label{leray}
	    |u(t)|_0^2+2\nu\int_{t'}^t\|u(s)\|_0^2ds\leq |u(t')|_0^2+2\int_{t'}^t(f(s),u(s))_0ds.
	\end{equation}
	for all $t\in I$, with $t'>t.$
\end{itemize}
\end{definition}
\begin{definition}
    A weak solution is said to be a strong/regular solution if it also belongs to $L^{\infty}(0,T; V_0)\cap L^2(0,T;D(A_0))$.\\
\end{definition}

\noindent
We have the following previously established results on existence and uniqueness.
\begin{theorem}\label{welldef1}
        Let $f\in L^2(0,T;V_0')$. Then, there exists a weak solution $u$ of \eqref{3dnav1}, satisfying all the properties given in Definition \ref{weakdef}.
\end{theorem}
\begin{theorem}\label{welldef2}
        Let $v$ be a strong solution to \eqref{3dnav1}. Then, there doesn't exist any other weak solution $u$ of \eqref{3dnav1}.
\end{theorem}
\begin{remark}
    For 3-D Navier-Stokes equation, we have existence of weak solution(Theorem \ref{welldef1}), but not uniqueness. Moreover, we have uniqueness of strong solutions(Theorem \ref{welldef2}), but not existence. 
\end{remark}
\subsection{Weak Attractor}
Despite the lack of a well-posedness result for the three-dimensional Navier-
Stokes equations, it is still a natural question to ask what the dynamics and the asymptotic behaviors of their weak solutions are, despite the possibility that they are not unique with respect to the initial condition. In particular, it is natural to ask whether there exists some sort of global attractor in this case. Due to the lack of a well-defined semigroup associated with the solutions of the system, the classical theory of dynamical system does not apply directly. Nevertheless, it is still possible to adapt a number of results from the classical theory to this situation.\\
One of the first and main results in this direction was given in \cite{FT87}, in which an object called the weak global attractor was defined as follows.
\begin{definition}
    The Weak attractor for the Navier-Stokes equation, $\mathbb{A}$, denotes the set of $u_0\in H_0$ for which there exists a weak solution $u(t)$ of \eqref{3dnav1}, for $t\in \mathbb{R}$, such that 
    \begin{itemize}
        \item $u\in L^{\infty}(\mathbb{R};H_0)$
        \item $u(0)=u_0.$
    \end{itemize}
\begin{remark}
    The weak global attractor for the 3-D Navier-Stokes operator has the following properties:
    \begin{itemize}
        \item  For every weak solution $u$ of \eqref{3dnav1} on the time interval $(0,\infty)$, we have 
    \begin{equation*}
        u(t)\to \mathbb{A}\text{ weakly in }H_0, \text{ as }t\to \infty.
    \end{equation*}
     \item $\mathbb{A}$ is weakly compact in $H_0.$
     \item $\mathbb{A}$ is invariant in the sense that if $u_0\in \mathbb{A}$ and $u$ is a global weak solution uniformly bounded in $H_0$ with $u(0)=u_0$ then $u(t)\in\mathbb{A}$ for all $t\in\mathbb{R}.$
    \end{itemize}
\end{remark}
\end{definition}
\subsection{Determining Map}
We begin by introducing spaces that contain the domain and ranges of the determining map $W$ and by introducing the evolution equation which yields the definition of $W$.\\
We denote by $L^2_b(\mathbb{R};D(A_0))$ the functions in $L^2(\mathbb{R};D(A_0))$ which are translation bounded, i.e.,
\begin{equation}
   \sup_{s\in\mathbb{R}} \int_s^{s+1}|A_0(u(r))|_{0}^2dr<\infty.
\end{equation}
Similarly, $L^2_b(\mathbb{R};V)$ denotes the functions in $L^2(\mathbb{R};V)$ which satisfy $$\sup_{s\in\mathbb{R}}\int_s^{s+1}\|u(r)\|_{0}^2dr<\infty.$$
Let $$Y=C_b(\mathbb{R};V_0)\cap L^2_b(\mathbb{R};D(A_0))\text{ and }Z=C_b(\mathbb{R};H_0)\cap L^2_b(\mathbb{R};V_0).$$
where, $C_b(I;B)$ denotes the space of all bounded and continuous functions over the interval $I$ with values in the Banach space $B$.
$Y$ and $Z$ are Banach spaces with norms
\begin{equation}
    \|u\|_Y=\left\{\|u\|_0^2+\sup_{s\in\mathbb{R}}\int_s^{s+1}|A_0(u(r))|_{0}^2dr\right\}^{\frac{1}{2}}.
\end{equation}
and
\begin{equation}
    \|u\|_Z=\left\{|u|_0^2+\sup_{s\in\mathbb{R}}\int_s^{s+1}\|u(r)\|_{0}^2dr\right\}^{\frac{1}{2}}
\end{equation}
respectively.\\
For $\tau\in I\subset\mathbb{R}$ and a Banach space $B$, we define the time translation $\tau_{\sigma}:C(I;B)\to C(I;B)$ as
\begin{equation}
    \tau_{\sigma}(u(t))=u(t+\sigma).
\end{equation}
Moreover, let $X$ be the Banach space $$X=C_b(\mathbb{R};(\dot{H}^1(\Omega))^3)\cap V_0$$ equipped with the norm 
\begin{equation}
    \|v\|_X=\sup_{s\in\mathbb{R}}\|v(s)\|_0.
\end{equation}
The observed spatial coarse-mesh data is denoted by $v(t)$. For the purpose of data assimilation, we consider the case $v\in X$ with $\|v\|_X\leq\rho$ for some $\rho>0$. We use $B_X(\rho)$ to denote the closed ball in $X$ of radius $\rho$ centered at $0$.\\
The data assimilation algorithm is given by the solution $w$ of the equation
\begin{equation}\label{datnav}
\begin{split}
     \frac{\partial w}{\partial t}+B_0(w,w)+\nu A_0(w)&=f + \mu(v-I_h(w))\\
     \nabla\cdot w&=0
    \end{split}
\end{equation}
Observe that the above system is not an initial value problem but an evolution equation for all $t\in \mathbb{R}.$
The Galerkin approximation of \eqref{datnav} is obtained by applying $P_n$ and is given by
\begin{equation}\label{datnav1}
\begin{split}
     \frac{\partial w_n}{\partial t}+B_0(w_n,w_n)+\nu A_0(w_n)&=P_n f + \mu(P_n(v)-I_h(w_n))\\
      \nabla\cdot w_n&=0
    \end{split}
\end{equation}
where $u_n$ is as in \eqref{3dnav2}.\\

\begin{theorem}\label{attexist}
		Let $u$ be the solution to \eqref{main1}-\eqref{main2} on the weak attractor for $t\in \mathbb{R}$, $v\in B_X(\rho)$ for some $\rho>0$ and $I_h$ be any type 1 interpolant. Define
		\begin{equation}
		    M_h^2=8\left(\frac{|f|^2}{\nu^2\lambda_1}+\rho^2\right).
		\end{equation}
		Also, let $h_0>0$ be defined as 
		\begin{equation}
		    h_0^2=\frac{1}{4c\lambda_1}.
		\end{equation}
		Assume that for some $h\leq h_0$,
		\begin{equation}
		\frac{2cM_h^4}{\nu^3}\leq \frac{\nu}{4ch^2}.
		\end{equation}
		Let $\mu$ be chosen such that
		\begin{equation}\label{attcondn}
		\max\left\{\frac{2cM_h^4}{\nu^3},\frac{\nu}{4ch_0^2}\right\}\leq\mu\leq \frac{\nu}{4ch^2}.
		\end{equation}
		Then there exists a unique global solution $(w,\eta)$ of \eqref{main4}-\eqref{main5} such that
		\begin{equation}\label{strongatt}
		    w\in L^{\infty}(\mathbb{R};V_0)\cap L^{2}(\mathbb{R};D(A_0)).
		\end{equation}
		Moreover, the following bounds hold 
		\begin{enumerate}
		    \item $\|w\|_0\leq M_h.$\label{one}\\
		    \item $\displaystyle\int_t^{t+1} |A_0w(t)|_0^2dt\leq \frac{4}{\nu}\left(\frac{1}{\nu}|f|_0^2+\mu\rho^2\right)$\label{two}
		\end{enumerate}
		Also, consider $v_1,v_2\in B_X(\rho)$ and let $w_1$ and $w_2$ be solutions to \eqref{datnav} corresponding to inputs $v_1$ and $v_2$ respectively. Denote $\tilde{w}=w_1-w_2$ and $\tilde{v}=v_1-v_2.$ Then
		
		\begin{enumerate}[resume]
		  \item $| \tilde{w}(t)|_0^2\leq 4\|\tilde{v}\|_X^2$\label{three}\\
		  \item $\displaystyle\int_s^{s+1}\|\tilde{w}(t)\|^2_0dt\leq \frac{4\mu}{\nu}\|\tilde{v}\|_X^2$\label{four}
		\end{enumerate}
	\end{theorem}
	\begin{proof}
	    The proof existence of global solution satisfying \eqref{strongatt} is obtained by showing that the solution, $w_n$ to the Galerkin system \eqref{datnav1} with $w(-N)=0$ satisfies $\eqref{one}$ and $\eqref{two}$ on the time interval $[-N,\infty)$(the proof for which follows exactly that of Theorem \ref{velreg1}) and  extracting a subsequence via the diagonal process and then passing to the limit. We proceed as in Theorem \ref{velreg1}. We will first look at equation \eqref{datnav1} on the time interval $[-N,\infty]$ with the initial condition $w_n(-N)=0.$\\
	    
	    \noindent
	    Taking the inner product of \eqref{datnav1} with $A_0w_n$, we obtain
	    \begin{equation}\label{navest}
	    \begin{split}
	    \frac{1}{2}\frac{d}{dt}\|w_n\|_0^2 +\nu|A_0w_n|_0^2=&-(B_0(w_n,w_n),A_0(w_n))_0+\mu(w_n-I_h(w_n),A_0(w_n))_0\\
	    &-\mu\|w_n\|_0^2+\mu(P_n(v),A_0(w_n))_0+(f,A_0(w_n))
	    \end{split}
	    \end{equation}
	    We bound each term below.
	    First, applying \eqref{nolinest1} and Young’s inequality, we have
	    \begin{align*}
	    |(B_0(w_n,w_n),A_0(w_n))|_0&\leq c\|w\|_0^{3/2}|A_0(w_n)|_0^{3/2}\\
        &\leq \frac{c}{\nu^3}\|w_n\|_0^6+\frac{\nu}{4}|A_0(w_n)|_0^2.
    	\end{align*}
	    Next, from \eqref{modalest}, Cauchy-Schwartz, Young’s inequality and the second inequality in \eqref{vstrong2}, we have
    	\begin{align*}
	    |\mu(w_n-I_h(w_n),A_0(w_n))|_0&\leq\mu ch\|w_n\|_0|A_0(w_n)|_0\\
    	&\leq \frac{\mu^2ch^2}{\nu}\|w_n\|_0^2+\frac{\nu}{4}|A_0(w_n)|_0^2\\
     	&\leq \frac{\mu}{4}\|w_n\|_0^2+\frac{\nu}{4}|A_0(w_n)|_0^2,
    	\end{align*}
    	Applying \eqref{A} and Cauchy-Schwarz, we obtain
    	\begin{align*}
    	\mu(P_n(v),A_0(w_n))_0&\leq \mu\|P_n(v)\|_0^2+\frac{\mu}{4}\|w_n\|_0^2.
    	\end{align*}
    	Lastly, applying Cauchy-Schwartz and Young’s inequality, we obtain
    	\begin{equation*}
    	    |(f,A_0w_n)_0|\leq |f|_0|A_0w_n|_0\leq\frac{1}{\nu}|f|_0^2+\frac{\nu}{4}|A_0w_n|^2
    	\end{equation*}
    	Inserting the above estimate into \eqref{navest}, we obtain
    	\begin{equation}\label{contra}
    	    \frac{d}{dt}\|w_n\|_0^2+\left(\mu-\frac{c}{\nu^3}\|w_n\|_0^4\right)\|w_n\|_0^2+\frac{\nu}{2}|A_0w_n|_0^2\leq \frac{2}{\nu}|f|_0^2+2\mu\|P_n(v)\|_0^2
    	\end{equation}
    	Let $[-N, T_1]$ be the maximal interval on which $\|w_n(t)\|\leq M_{h}$ for $t\in[-N, T_1]$, 
	where $M_{h}$ as in \eqref{attcondn}. Note that $T_1 > -N$ exists because we have $w_n(-N)=0$. Assume
	that $T_1<\infty$ . Then by continuity, we must have $\|w_n(T_1)\|=M_{h}$. Applying \eqref{unibound},\eqref{weakbound}, the first inequality in \eqref{vstrong2} and dropping all terms except the first and the last term on the LHS of \eqref{contra}, we obtain
	\begin{equation}\label{contra2}
	\frac{d}{dt}\|w_n\|_0^2 +\frac{\mu}{2}\|w_n\|_0^2\leq \frac{2}{\nu}|f|_0^2+2\mu\sup_{t\in\mathbb{R}}\|P_n(v)\|_0^2.
	\end{equation}
		Applying Gronwall, we obtain, for all $t\in [-N,T_1],$
	\begin{align*}
	\|w_n\|_0^2 \leq 4\left(\frac{|f|^2}{\nu^2\lambda_1}+\|(v)\|_X^2\right)\left(1-e^{-\frac{\mu}{8}(n+t)}\right)&\leq 4\left(\frac{|f|^2}{\nu^2\lambda_1}+\rho^2\right)\\
	&\leq \frac{1}{2}M_{h}^{2}.
	\end{align*}
	This contradicts the fact that $\|w_n(T_1)\|=M_{h}$. Therefore $\|w_n(t)\|\leq M_{h}$ for all $t\in[-n, \infty).$
    Dropping all terms except the last term on the LHS of \eqref{contra} and integrating both sides over the interval $[t, t+1]$, it follows that
    	\begin{equation}
    	   \int_t^{t+1} |A_0w_n|_0^2\leq \frac{4}{\nu}\left(\frac{1}{\nu}|f|_0^2+\mu\rho^2\right)
    	\end{equation}
    Therefore, we have a sequence of solutions $w_n$ that satisfy $\eqref{one}$ and $\eqref{two}$ on the interval $[-N,\infty)$. One can now extract a convergent subsequence via the diagonal process and pass through the limit to show that the limit $w$ is a solution to \eqref{datnav} and satisfies \eqref{strongatt}, $\eqref{one}$ and $\eqref{two}$.\\
    
    \noindent
    We now proceed to prove \eqref{three} and \eqref{four}. Let $w_{1,n}$ and $w_{2,n}$ satisfy \eqref{datnav1} with $v=v_1$ and $v=v_2$ respectively. Then $\tilde{w}_n=w_{1,n}-w_{2,n}$ and $\tilde{v}_n=P_n(\tilde{v})$ satisfy 
   	\begin{equation}\label{uniq}
	\begin{split}
	\frac{d\tilde{w}_n}{dt}+\nu A_0(\tilde{w}_n)+B_0(\tilde{w}_n,w_{1,n})+B_0(w_{2,n},\tilde{w}_n)&=\mu P_n(\tilde{v}-I_h(\tilde{w}))\\
	&=\mu \tilde{v}_n+\mu(\tilde{w}_n-I_h(\tilde{w}_n))\\
	&-\mu\tilde{w}_n
	\end{split}
	\end{equation}
	Taking the inner product of \eqref{uniq} with $\tilde{w}$, we obtain
	\begin{equation}\label{mmmbop}
	\begin{split}
	\frac{1}{2}\frac{d}{dt}|\tilde{w}_n|_0^2+\nu\|\tilde{w}_n\|_0^2+\mu|\tilde{w}_n|_0^2&\leq |(B_0(\tilde{w}_n,w_{1,n}),\bar{w}_n|)|_0\\
	&+|\mu(\tilde{w}_n-I_h(\tilde{w}_n))|_0|\bar{w}_n|_0\\
	&+\mu|\tilde{v}_n|_0|\bar{w}_n|_0.
	\end{split}
	\end{equation}
	We bound each term on the RHS.\\
	First, applying \eqref{nolinest2}, Cauchy-Schwartz and Young’s inequality, we obtain
	\begin{equation}
	\begin{split}
	|(B_0(\tilde{w}_n,w_{1,n}),\tilde{w}_n|)|_0&\leq c|\tilde{w}_n|_0^{1/2}\|\tilde{w}_n\|_0^{3/2}\|w_{1,n}\|_0\\
	&\leq \frac{c}{\nu^3}\|w_{1,n}\|_0^4|\tilde{w}_n|_0^2+\frac{\nu}{2}\|\tilde{w}_n\|_0^2.
	\end{split}
	\end{equation}
	Applying \eqref{intest}, Cauchy-Schwartz, Young’s inequality and \eqref{attcondn}, we obtain
	\begin{equation}
	\begin{split}
	|\mu(\tilde{w}_n-I_h(\tilde{w}_n))|_0|\tilde{w}_n|_0&\leq \mu ch\|\tilde{w}_n\|_0|\tilde{w}_n|_0\\
	&\leq \mu ch^2\|\tilde{w}_n\|_0^2+\frac{\mu}{4}|\tilde{w}_n|_0^2\\
	&\leq \frac{\nu}{4}\|\tilde{w}_n\|_0^2+\frac{\mu}{4}|\tilde{w}_n|_0^2.
	\end{split}
	\end{equation}
	Lastly, applying Young’s inequality we have
	\begin{equation}
	\begin{split}
	\mu|\tilde{v}_n|_0|\bar{w}_n|_0\leq \mu|\tilde{v}_n|_0^2+\frac{\mu}{4}|\bar{w}_n|_0^2,
	\end{split}
	\end{equation}
	Inserting the above estimates into\eqref{mmmbop}, we obtain,
	\begin{equation}\label{mmmbop2}
	\begin{split}
	\frac{d}{dt}|\tilde{w}_n|_0^2+\frac{\nu}{2}\|\tilde{w}_n\|_0^2+\left(\mu-\frac{c}{\nu^3}\|w_{1,n}\|_0^4|\right)|\tilde{w}_n|_0^2&\leq 2\mu|\tilde{v}_n|_0.
	\end{split}
	\end{equation}
	From proof of \eqref{one}, we know that $\|w_{1,n}(t)\|_0\leq M_h~~\forall t\in\mathbb{R}$. Hence, applying \eqref{attcondn} and integrating on the interval $[\sigma,t]$, we obtain
	\begin{equation}
	\begin{split}
	|\tilde{w}_n(t)|_0^2&\leq |\tilde{w}_n(\sigma)|^2e^{\frac{\mu}{2}(\sigma-t)}+4\|\tilde{v}_n\|_X^2\left(1-e^{\frac{\mu}{2}(\sigma-t)}\right).
	\end{split}
	\end{equation}
	Letting $\sigma\to-\infty$, we obtain
	\begin{equation}\label{fi}
	   | \tilde{w}_n(t)|_0^2\leq 4\|\tilde{v}_n\|_X^2
	\end{equation}
	Keeping only the second term on the LHS of \eqref{mmmbop2} and integrating on the interval$[s,s+1],$ we obtain
	\begin{equation}\label{se}
	    \int_s^{s+1}\|\tilde{w}_n(t)\|^2_0dt\leq \frac{4\mu}{\nu}\|\tilde{v}_n\|_X^2
	\end{equation}
	Taking the limit $n\to\infty$ in \eqref{fi} and \eqref{se}, we obtain \eqref{three} and \eqref{four}.
    \end{proof}
    
    \begin{definition}\label{mappy}
	Consider $\rho>0$, $\mu>0$ and $h>0$ satisfying the hypothesis of Theorem \ref{attexist}. Then, the determining map $W:B_X(\rho)\to Y$ is given  by $W(v)=w.$
	\end{definition}
	\begin{remark}
	    Observe that $Y\subset Z.$ Hence, from \eqref{three} and \eqref{four} in Theorem \ref{attexist}, we may conclude that $W:B_X(\rho)\to Z$ is Lipschitz continuous. Lipschitz continuity in turn implies uniqueness of solutions.
	\end{remark}
	\begin{corollary}
	    The determining map $W$ defined in Definition \ref{mappy}, in addition to being Lipschitz continuous has the following properties
	    \begin{enumerate}
	        \item Let $v_1,v_2\in B_{X}(\rho)$, Then, $W(v_1)=W(v_2)$ iff $P_{\sigma}(v_1-v_2)=P_{\sigma}(\bar{v})=0.$\label{55}\\
	    \item For every $\sigma\in \mathbb{R}$,
	    \begin{equation}
	        W\circ\tau_\sigma(v)=\tau_\sigma\circ W(v).
	    \end{equation} \label{666}
	    \end{enumerate}
	\end{corollary}
	\begin{proof}
	    To prove \eqref{55}, we note that $\bar{w}=w_1-w_2=0$, where $w_1=W(v_1)$ and $w_2=W(v_2)$. Since $w_1$ and $w_2$ are regular, the term $\bar{w}$ is differentiable a.e on $\mathbb{R}$. From \eqref{datnav}, we may write
	\begin{equation}
	    	\frac{\partial\bar{w}}{\partial t}+\nu A_0(\bar{w})+B_0(\bar{w},w_1)+B_0(w_2,\bar{w})=\mu P_\sigma(\bar{v}-I_h(\bar{w}))
	\end{equation}
	Letting $\bar{w}=0$, we obtain $P_\sigma(\bar{v})=0.$ If $P_\sigma(\bar{v})=0$, We obtain $W(v_1)=W(v_2)$ from the Lipschitz continuity of $W.$\\
	To show \eqref{666}, we observe that $\tau_\sigma\circ W(v)$ is a solution of \eqref{datnav} corresponding to $\tau_\sigma(v)$. From the Lipschitz property of the map $W$, we have uniqueness of solution. Hence $$ W\circ\tau_\sigma(v)=\tau_\sigma\circ W(v).$$
	\end{proof}
	\subsection{New Regularity Criterion On The Weak Attractor}
    In this section we establishes that the Leray-Hopf weak solution $u$ to \eqref{3dnav1} on weak attractor is in fact a strong solution. We will use the following four lemmas to prove Theorem \ref{yawza}, which will help us obtain our regularity criterion(given by Theorem \ref{criterion}).\\
    The first lemma is obtained by repeating the arguments of Lemma \ref{sather-serrin}.
	\begin{lemma}\label{sather-serrin2}
	    Let $u$ be the general Leray-Hopf weak solution of \eqref{3dnav} and $w$ be the strong solution to \eqref{datnav} as given in  \cite{Nav3d}. Let $\tilde{w}=w-u$ . Then, 
	    \begin{equation}
	    \begin{split}
	         |\tilde{w}(t)|_0^2+2\nu\int_0^t\|\tilde{w}(s)\|_0^2ds&\leq |\tilde{w}(0)|_0^2+2\int_0^tb(\tilde{w}(s),\tilde{w}(s),w(s))ds\\
            &-2\mu\int_0^t\left(I_h(\tilde{w}(s),\tilde{w}(s))\right)_0ds	  
	    \end{split}
        \end{equation}
	\end{lemma}

	\begin{lemma}\label{agree}
	    Let $u$ be a general Leray-Hopf weak solution of \eqref{3dnav} and $w$ be a strong solution to \eqref{datnav} as given in \cite{Nav3d}. Also, let $\tilde{w}=w-u.$ Then,
	    \begin{equation}
	        |\tilde{w}(t)|_0^2\leq 4M^2e^{-\mu\left(t-\frac{1}{\mu}\right)}.
	    \end{equation}
	    where 
	    \begin{equation}\label{M}
	        M=\max\left\{\displaystyle\sup_{0\leq t\leq \frac{1}{\mu}}|w(t)|_0,\sup_{0\leq t\leq \frac{1}{\mu}}|u(t)|_0\right\}
	    \end{equation}.
	\end{lemma}
	\begin{proof}
	    Repeating the arguments of Lemma \ref{sather-serrin}, we have
	    \begin{equation}\label{oho}
	    \begin{split}
	        |\tilde{w}(t)|_0^2+2\nu\int_0^t\|\tilde{w}(s)\|_0^2ds&\leq|\tilde{w}(0)|_0^2+ 2\int_0^t (B_0(\tilde{w}(s),\tilde{w}(s)),w(s))_0ds
            -2\mu\int_0^t\left(I_h(\tilde{w}(s),\tilde{w}(s))\right)_0ds\\
            &\leq|\tilde{w}(0)|_0^2+ 2\int_0^t (B_0(\tilde{w}(s),\tilde{w}(s)),w(s))_0ds\\
            &+2\mu\left(\int_0^t\left(\tilde{w}-I_h(\tilde{w}(s),\tilde{w}(s))\right)_0ds\right)-2\mu\int_0^t|\tilde{w}|_0^2ds
            \end{split}
        \end{equation}
        We bound each term on the RHS.
         First, applying \eqref{nolinest2}, Cauchy-Schwartz and Young’s inequality, we obtain
	\begin{equation}
	\begin{split}
	|(B_0(\tilde{w},\tilde{w}),w)_0|=(B_0(\tilde{w},w),\tilde{w}|)_0|
	&\leq c|\tilde{w}|_0^{1/2}\|\tilde{w}\|_0^{3/2}\|w\|_0\\
	&\leq \frac{c}{\nu^3}\|w\|_0^4|\tilde{w}|_0^2+\frac{\nu}{2}\|\tilde{w}\|_0^2\\
	&\leq \frac{c}{\nu^3}M_{h,u}^4|\tilde{w}|_0^2+\frac{\nu}{2}\|\tilde{w}\|_0^2
	\end{split}
	\end{equation}
	    Applying \eqref{intest}, Cauchy-Schwartz, Young’s inequality and \eqref{vstrong2}, we obtain
	\begin{equation}
	\begin{split}
	2|\mu(\tilde{w}-I_h(\tilde{w}))|_0|\tilde{w}|_0&\leq 2\mu ch\|\tilde{w}\|_0|\tilde{w}|_0\\
	&\leq 2\mu ch^2\|\tilde{w}\|_0^2+\frac{\mu}{2}|\tilde{w}|_0^2\\
	&\leq \frac{\nu}{2}\|\tilde{w}\|_0^2+\frac{\mu}{2}|\tilde{w}|_0^2.
	\end{split}
	\end{equation}
	$w$, as shown in \cite{Nav3d} satisfies the hypothesis of Theorem \ref{attexist}.
    Inserting the estimates into \eqref{oho} and applying \eqref{attcondn},  we see that
    \begin{equation}\label{synchro}
         |\tilde{w}(t)|_0^2+\nu\int_0^t\|\tilde{w}(s)\|_0^2ds+\mu\int_0^t|\tilde{w}|_0^2ds\leq |\tilde{w}(0)|_0^2
    \end{equation}
	Dropping the second term on the LHS and changing the interval of integration from $[0,t]$ to $[t-1/\mu,t]$, we have
	 \begin{equation}
         |\tilde{w}(t)|_0^2+\mu\int_{t-\frac{1}{\mu}}^t|\tilde{w}|_0^2ds\leq |\tilde{w}(t-1/\mu)|_0^2
    \end{equation}
	Now, applying Lemma \ref{titi}, we obtain
	\begin{equation}
	\begin{split}
	     |\tilde{w}(t)|_0^2&\leq \mu e^{-\mu\left(t-\frac{1}{\mu}\right)}\int_{0}^{\frac{1}{\mu}}|\tilde{w}(\tau)|_0^2d\tau\\
	     &\leq 4M^2e^{-\mu\left(t-\frac{1}{\mu}\right)}
	\end{split}
	\end{equation}
	\end{proof}
	\begin{lemma}\label{tau}
	    Let $u$ be a Leray-Hopf weak solution of \eqref{3dnav} on the weak attractor and $w$ be the strong solution of \eqref{datnav} as given in the Theorem \ref{attexist}. Let $\tilde{w}=w-u$ and $\displaystyle\tilde{w}_{\tau}(t)=\tilde{w}\left(t-\tau-1/\mu\right)$. Then
	    \begin{equation}
	          |\tilde{w}_{\tau}(t)|_0^2\leq 4M^2e^{-\mu\left(t-\frac{1}{\mu}\right)}.
	    \end{equation}
	    where $M$ is as in \eqref{M}.
	\end{lemma}
	\begin{proof}
	    Dropping the middle term on the LHS of \eqref{synchro}, we obtain
	    \begin{equation}
         |\tilde{w}(t)|_0^2+\mu\int_0^t|\tilde{w}(s)|_0^2ds\leq |\tilde{w}(0)|_0^2
    \end{equation}
    Changing the interval of integration from $[0,t]$ to $[t-\tau-2/\mu, t-\tau-1/\mu]$, we obtain
     \begin{equation}\label{synchro1}
         |\tilde{w}( t-\tau-1/\mu)|_0^2+\mu\int_{ t-\tau-2/\mu}^{ t-\tau-1/\mu}|\tilde{w}(s)|_0^2ds\leq |\tilde{w}( t-\tau-2/\mu)|_0^2
    \end{equation}
   Applying the definition of $\tilde{w}_{\tau}$, \eqref{synchro1} can be rewritten as
    \begin{equation}
         |\tilde{w}_{\tau}(t)|_0^2+\mu\int_{t-1/\mu}^t|\tilde{w}_{\tau}(s)|_0^2ds\leq |\tilde{w}_{\tau}(t-1/\mu)|_0^2.
    \end{equation}
    Proceeding exactly as in Lemma \ref{agree}, we obtain the statement of the lemma.
	\end{proof}
    \begin{theorem}\label{yawza}
        Let $u$ be a general Leray-Hopf weak solution of \eqref{3dnav} and $w$ be a strong solution to \eqref{datnav} as given in Theorem \ref{attexist}. Also, let $\tilde{w}=w-u$ and $\displaystyle\tilde{w}_{\tau}(t)=\tilde{w}\left(t-\tau-1/\mu\right)$. Then, $|\tilde{w}(t)|_0=0~\forall t\in \mathbb{R}.$
    \end{theorem}
    \begin{proof}
         It is enough to show that for any $t_0\in \mathbb{R}$ ,
         $|\tilde w(t_0)|_0^2 \leq \epsilon$. Given $t_0\in \mathbb{R}$, let $\tau\in \mathbb{R}$ be such that 
         \begin{equation}\label{yaha}
            \left(4M^2e^{-\mu\tau}\right)e^{-\mu\left(t_0\right)}\leq \epsilon,
         \end{equation}
         where $M$ is as in \eqref{M}. Then, from Lemma \ref{tau}, we have
         \begin{equation}
	          |\tilde{w}_{\tau}(t_0+\tau+1/\mu)|_0^2\leq \left(4M^2e^{-\mu\tau}\right)e^{-\mu\left(t_0\right)}.
	    \end{equation}
	    Applying \eqref{yaha} and the definition of $\tilde{w}_{\tau}$, we obtain
	    \begin{equation}
	        |\tilde{w}(t_0)|_0^2\leq\epsilon.
	    \end{equation}
    \end{proof}
    \begin{remark}
        From the above theorem we see that on the weak attractor, $u=w$, where $u$ is a general Leray-Hopf weak solution to \eqref{3dnav} and $w$ be a strong solution to \eqref{datnav}. This would in turn mean that $u$ on the weak attractor can be shown to be regular whenever we can can construct such a $w$, giving us a regularity criterion summarized in the thorem below.
    \end{remark}
    \begin{theorem}\label{criterion}
         Let $u$ be a  Leray-Hopf weak solution of \eqref{3dnav1} on the weak attractor. Let $h_0>0$ be defined as
         \[
             h_0^{-2}=\max\left\{\frac{1}{4c\lambda_1},\frac{32c|f|^4}{\nu^8\lambda_1^2}\right\}\quad \text{($f$ is the body force)}.
         \]
         Assume there exists $0<h\leq h_0$ for which $$\frac{2cM_{h}^4}{\nu^3}\leq \frac{\nu}{16ch^2}.$$ Then $u(t)$ is regular for all $t\in\mathbb{R}$.
    \end{theorem}
    
\section{Weakened Regularity Condition}\label{weakenn}
In Theorem \ref{attexist}, the condition we imposed on the data to show regularity of $w$ can be weakened. In this section, we will show that the weakened condition on the data is sufficient to obtain the same regularity result as that in Theorem \ref{attexist}\eqref{one}.\\
We define $K_h:=K_{h,\tau_0}$, for $\tau_0>0$ and Holder conjugates $p\geq 3$ and $q\geq 1$, as
	\begin{equation}\label{k}
	    K_h=\sup_{t\in[-\infty,T]}\left(\int_t^{t+\tau_0}\|I_hu\|^{2p} ds\right)^{1/2p}.
	\end{equation}
The weakened condition is obtained in terms of $K_h$ instead of $\rho$.
\begin{theorem}
	Let $\tilde{I}_h$ be a modified general type-\rom{1}  interpolant, $K_h$ be defined as in \eqref{k} and $0\leq T\leq\infty$. 
	Let $h_0>0$ be given by
	\begin{equation}  
        h_0^{-2}=\max\left\{4c\lambda_1,\frac{1024c|f|^4}{\nu^8\lambda^2}\right\}
	\end{equation}
	Assume there exists $\tau_0>0$ and $0<h\leq h_0$ such that $K_h<\infty$ and
	\begin{equation}
	    \left(\frac{32CK_h^4}{q^{\frac{2}{q}}}\right)^{\frac{p}{p-2}}\leq\frac{\nu}{4ch^2}.
	\end{equation}
	Let $\mu$ be chosen such that
	\begin{equation}\label{bigcond}
	    \max\left\{\frac{\nu}{4ch_0^2},\left(\frac{32CK_h^4}{q^{\frac{2}{q}}}\right)^{\frac{p}{p-2}}\right\}\leq \mu\leq\frac{\nu}{4ch^2}.
	\end{equation}
	Then, the data assimilated fluid velocity is regular and satisfies \begin{equation}
	    \|w\|_0\leq M_{h},
	\end{equation}
	where 
	\begin{equation}\label{Mh}
	    M_{h}^2=\frac{8|f|^2}{\lambda_1\nu^2}+\frac{2CK_h^2\mu^{1/p}}{q^{1/q}}\left(\frac{2}{1-e^{-\frac{\nu\lambda_1 p}{4} \tau_0}}\right)^{1/p}.
	\end{equation}
\end{theorem}

\begin{proof}
Taking the inner product of the 3-D Navier-Stokes equation with $A_0w_n$, we obtain
	    \begin{equation}\label{navest1}
	    \begin{split}
	    \frac{1}{2}\frac{d}{dt}\|w_n\|_0^2 +\nu|A_0w_n|_0^2=&-(B_0(w_n,w_n),A_0(w_n))_0+\mu(w_n-I_h(w_n),A_0(w_n))_0\\
	    &-\mu\|w_n\|_0^2+\mu(I_hu,A_0(w_n))_0+(f,A_0(w_n))
	    \end{split}
	    \end{equation}
	    We bound each term below.
	    
	    \begin{align*}
	    |(B_0(w_n,w_n),A_0(w_n))|_0&\leq c\|w\|_0^{3/2}|A_0(w_n)|_0^{3/2}\\
        &\leq \frac{c}{\nu^3}\|w_n\|_0^6+\frac{\nu}{4}|A_0(w_n)|_0^2.
    	\end{align*}
	    
    	\begin{align*}
	    |\mu(w_n-I_h w_n,A_0(w_n))|_0&\leq\mu ch\|w_n\|_0|A_0(w_n)|_0\\
    	&\leq \frac{\mu^2ch^2}{\nu}\|w_n\|_0^2+\frac{\nu}{4}|A_0(w_n)|_0^2\\
     	&\leq \frac{\mu}{4}\|w_n\|_0^2+\frac{\nu}{4}|A_0(w_n)|_0^2,
    	\end{align*}
    	
    	\begin{align*}
    	\mu(I_hu,A_0(w_n))_0&\leq \mu\|I_hu\|_0^2+\frac{\mu}{4}\|w_n\|_0^2.
    	\end{align*}
    	Lastly, applying Cauchy-Schwartz and Young’s inequality, we obtain
    	\begin{equation*}
    	    |(f,A_0w_n)_0|\leq |f|_0|A_0w_n|_0\leq\frac{1}{\nu}|f|_0^2+\frac{\nu}{4}|A_0w_n|^2
    	\end{equation*}
    	Inserting the above estimate into \eqref{navest1}, we obtain
    	\begin{equation}\label{contraa}
    	    \frac{d}{dt}\|w_n\|_0^2+\left(\mu-\frac{c}{\nu^3}\|w_n\|_0^4\right)\|w_n\|_0^2+\frac{\nu}{2}|A_0w_n|_0^2\leq \frac{2}{\nu}|f|_0^2+2\mu\|I_hu\|_0^2
    	\end{equation}
    	Let $[0, T_1]$ be the maximal interval on which $\|w_n(t)\|\leq M_{h}$ for $t\in[0, T_1]$. Note that $T_1 > 0$ exists because we have $w_n(0)=0$. Assume
	that $T_1<T$ . Then by continuity, we must have $\|w_n(T_1)\|=M_{h}$. Applying   \eqref{bigcond}and \eqref{Mh} to \eqref{contraa} and dropping the last term on the LHS of \eqref{contraa}, we obtain
	\begin{equation}\label{contraa2}
	\frac{d}{dt}\|w_n\|_0^2 +\frac{\mu}{2}\|w_n\|_0^2\leq \frac{2}{\nu}|f|_0^2+2\mu\|P_n(v)\|_0^2.
	\end{equation}
	Applying Gronwall's and Holder's inequalities, for any $t\in[0,T_1]$ we obtain
	\begin{equation}\label{wc}
	    \begin{split}
	        \|w_n(t)\|^2&\leq \frac{4|f|^2}{\mu\nu}+2\mu\int_0^t e^{-\frac{\mu}{2}(t-s)}\|I_hu\|^2ds\\
	        &\leq \frac{4|f|^2}{\mu\nu}+2\mu\left(\int_0^te^{-\frac{\mu q}{4} (t-s)}ds\right)^{1/q}\left(\int_0^t e^{-\frac{\mu p}{4} (t-s)}\|I_hu\|^{2p}ds\right)^{1/p}\\
	        &\leq  \frac{4|f|^2}{\mu\nu}+\frac{C\mu^{1/p}}{q^{1/q}}\left(\int_0^t e^{-\frac{\mu p}{4} (t-s)}\|I_hu\|^{2p}ds\right)^{1/p}
	    \end{split}
	\end{equation}
	Where $1\leq p,q\leq\infty$ are Holder conjugates.
	We now try to bound the second term on the RHS of the above inequality.\\
	
	\noindent
	\textit{\textbf{Case 1}}: $t\geq \tau_0$.\\
	
	\noindent
	Let $\displaystyle k=\floor*{\frac{t}{\tau_0}}$ be the largest integer such that $k\tau_0\leq t.$ Therefore,
	\begin{equation}\label{k1}
	    t=k\tau_0+\epsilon,~~\text{where}~~ 0\leq\epsilon\leq\tau_0.
	\end{equation}
	Let $\displaystyle \alpha=\frac{p\mu}{4}$. Then, from \eqref{k} and \eqref{k1}, we may write
	\begin{equation}\label{case1}
	    \begin{split}
	        \int_0^t e^{-\frac{\mu p}{4} (t-s)}\|I_hu\|^{2p}ds&\leq \sum_{j=1}^k\int_{(j-1)\tau_0}^{j\tau_0}e^{-\alpha(t-s)}\|I_h(u)\|^{2p}ds+\int_{n\tau_0}^{t} e^{-\alpha(t-s)}\|I_h(u)\|^{2p}ds\\
	        &\leq e^{-\alpha\epsilon}\sum_{j=1}^k e^{-\alpha(n-j)\tau}\int_{(j-1)\tau_0}^{j\tau_0}\|I_h(u)\|^{2p}ds+\int_{n\tau_0}^{n\tau_0+\epsilon}\|I_h(u)\|^{2p}ds\\
	        &\leq e^{-\alpha\epsilon}K_h^{2p}\sum_{m=0}^{k-1}e^{m\alpha\tau_0}+K_h^{2p}.\\
	        &\leq \frac{2K_h^{2p}}{1-e^{\alpha\tau_0}}
	    \end{split}
	\end{equation}
	\textit{\textbf{Case 2}}: $t\leq \tau_0$.\\
	
	\noindent
	From\eqref{k}, we see that
	\begin{equation}\label{case2}
	    \begin{split}
	        \int_0^t e^{-\alpha (t-s)}\|I_hu\|^2ds&\leq \int_0^t \|I_hu\|^2ds\\
	        &\leq \int_0^{\tau_0} \|I_hu\|^2ds\\
	        &\leq K_h^{2p}\\
	        &\leq \frac{2K_h^{2p}}{1-e^{\alpha\tau_0}}
	    \end{split}
	\end{equation}
	Applying \eqref{case1}, \eqref{case2} and \eqref{Mh} to \eqref{wc}, we obtain that for any $t\in[0,T_1]$
	\begin{equation}
	    \begin{split}
	        \|w_n(t)\|^2&\leq  \frac{4|f|^2}{\mu\nu}+\frac{CK_h^2\mu^{1/p}}{q^{1/q}}\left(\frac{2}{1-e^{-\frac{\mu p}{4} \tau_0}}\right)^{1/p}\\
	        &\leq \frac{4|f|^2}{\lambda_1\nu^2}+\frac{CK_h^2\mu^{1/p}}{q^{1/q}}\left(\frac{2}{1-e^{-\frac{\nu\lambda_1 p}{4} \tau_0}}\right)^{1/p}\\
	        &\leq\frac{1}{2}M_h^2.
	    \end{split}
	\end{equation}
	This contradicts the fact that $\|w_n(T_1)\|=M_{h}$. Therefore $T_1\geq T$ and consequently, $\|w_n(t)\|\leq M_{h}$ for all $t\in[0, T ].$ Passing to the limit as $n\to\infty$, we obtain the desired conclusion for $w$.
\end{proof}
Hence, similar to Theorem \ref{criterion}, using Lemma \ref{sather-serrin2} - Lemma \ref{tau}, we obtain the following regularity criterion:
    \begin{theorem}\label{criterion1}
         Let $u$ be a  Leray-Hopf weak solution of \eqref{3dnav1} on the weak attractor and $K_h$ be defined as in \eqref{k}. Let $h_0>0$ be defined as
         \[
              h_0^{-2}=\max\left\{4c\lambda_1,\frac{1024c|f|^4}{\nu^5\lambda^2}\right\}.
         \]
         Assume there exists $0<h\leq h_0$ for which $$\max\left\{\frac{\nu}{4ch_0^2},\left(\frac{32CK_h^4}{q^{\frac{2}{q}}}\right)^{\frac{p}{p-2}}\right\}\leq \mu\leq\frac{\nu}{4ch^2}.$$ Then $u(t)$ is regular for all $t\in\mathbb{R}$.
    \end{theorem}
    \begin{remark}
        Note that in Theorem \ref{criterion1}, the definition of $M_h$ is given by \eqref{Mh}, which is not the same as \eqref{mh}. Also, the regularity criterion given by Theorem \ref{criterion1} is in the spirit of the criterion given by Corollary 5.2 in \cite{TitiHolts}. However, our condition solely depends on the observed data.
    \end{remark}
\section{Appendix}
Recall that for volume interpolation, we divided our domain into smaller sub domains(cuboids $Q_{\alpha}$) of diameter $h$ and  indexed by the set $\mathcal{J}$. Let us define the set $\mathcal{E}\subset\mathcal{J}$ as 
\begin{equation}
    \mathcal{E}=\left\{\alpha=(\alpha_1,\alpha_2,\alpha_3)\in\mathcal{J}:\alpha_3=1\text{ or }\alpha_3=\sqrt[3]{N}\right\}.
\end{equation}
$\mathcal{E}$ represents the collection of sub domains touching the top and bottom boundaries($x_3=0$ and $x_3=1$). For sub domains $Q_{\alpha}$ with $\alpha\in\mathcal{E}$, we introduce a modification $Q_{\epsilon,\alpha}$ for $0<\epsilon<h$, given by
\begin{equation}
    Q_{\epsilon,\alpha}=\begin{cases}
\displaystyle \left((\alpha_1-1)*h_L,\alpha_1*h_L\right)\times\left((\alpha_2-1)*h_L,\alpha_2*h_L\right)\times(\epsilon,h_1)&,\text{ for }\alpha_3=1\\
\left((\alpha_1-1)*h_L,\alpha_1*h_L\right)\times\left((\alpha_2-1)*h_L,\alpha_2*h_L\right)\times(1-h_1,1-\epsilon)&,\text{ for }\alpha_3=\sqrt[3]{N},
\end{cases}
\end{equation}
where $h_x=x/\sqrt[3]{N}.$
We now define the $\tilde{I}_h$ as
\begin{equation}\label{modef}
\tilde{I}_h(v)(x) = \sum_{\alpha\in\mathcal{J}}\bar{v}_{\alpha}\phi_{\alpha}(x),~~ v\in H^1(\Omega) 
\end{equation}
where
\begin{equation}\label{phi}
\phi_{\alpha}=\rho_{\epsilon}*\psi_{Q_{\alpha}},~\rho_{\epsilon}(x)=\epsilon^{-3}\rho\left(x/{\epsilon}\right),~\bar{v}_{\alpha} = \frac{1}{|Q_\alpha|}\int_{Q_{\alpha}}v(x)dx.
\end{equation}

\begin{equation}
    \psi_{Q_{\alpha}}(x) =\begin{cases}
\displaystyle \chi_{Q_{\alpha}} &,\text{ for }\alpha\not\in\mathcal{E}\\
\chi_{Q_{\epsilon,\alpha}} &,\text{ for }\alpha\in\mathcal{E}
\end{cases}
\end{equation}

\begin{equation}\label{rho}
\rho(x) =\begin{cases}
\displaystyle K_{0}\exp\left(\frac{-1}{1-x^2}\right) &,\text{ for }|x|<1\\
~0 &,\text{ otherwise}
\end{cases}
\end{equation}
and 
\begin{equation*}
(K_0)^{-1}=\int_{|x|<1}\exp\left(\frac{-1}{1-x^2}\right)dx.
\end{equation*}
We set $\epsilon=h/10.$
\begin{remark}
    $\|\tilde{I}_h\|_0$ is well defined since the characteristic function of each sub domain has been mollified. The characteristic function of sub domains touching the top and bottom boundaries have been modified($Q_{\epsilon,\alpha}$) so that it's support is an "$\epsilon$ distance" away from the top and the bottom boundaries. This is done so that after mollification, the modified characteristic function respects the Dirichlet boundary condition at $x_3=0,1.$
    \end{remark}
\begin{lemma}\label{applemma}
	Let $\phi_{\alpha}$ and $\rho$ be defined as in \eqref{phi} and \eqref{rho} respectively. Then, for $i=1$, $2$ and $3$
	$$|\partial_{x_i}\phi_{\alpha}|^2\leq Ch\|\partial_{x_i}\rho\|_{L^\infty(\Omega)}^2$$
\end{lemma}
\begin{proof}
	Recall that $\phi_{\alpha}=\rho_{\epsilon}*\psi_{Q_{\alpha}}$. Applying Young's inequality for convolutions, we obtain
	\begin{equation}\label{convo}
	\begin{split}
	|\partial_{x_i}\left(\rho_{\epsilon}*\psi_{Q_\alpha}\right)|^2&=|\left(\partial_{x_i}\rho_{\epsilon}\right)*\psi_{Q_\alpha}|^2\\
	&\leq \left\|\partial_{x_i}\rho_{\epsilon}\right\|_{L^1(\Omega)}^2\left|\psi_{Q_\alpha}\right|^2
	\end{split}
	\end{equation}
	Now, we look at each term on the RHS of \eqref{convo}. Differentiating the second equation in \eqref{phi} with $x_i$ and using the fact that $\epsilon=h/10$, we obtain
	\begin{equation}\label{convo2}
	\begin{split}
	\left\|\partial_{x_i}\rho_{\epsilon}\right\|_{L^1(\Omega)}^2&=\left(\int_{|x|\leq\epsilon}\left|\epsilon^{-4}\partial_{x_i}\rho(x/\epsilon)\right|dx\right)^2\\
	&\leq C\epsilon^{-2}\|\partial_{x_i}\rho\|_{L^\infty(\Omega)}^2\\
	&\leq Ch^{-2}\|\partial_{x_i}\rho\|_{L^\infty(\Omega)}^2.
	\end{split}
	\end{equation}
	From the definition of $\chi_{Q_{\alpha}}$, we readily obtain
	\begin{equation}\label{convo3}
	\begin{split}
	\left|\psi_{Q_{\alpha}}\right|^2&\leq\int_{Q_{\alpha}}1^2dx=\frac{|\Omega|}{N}\leq h^3
	\end{split}
	\end{equation}
	Combining \eqref{convo}, \eqref{convo2} and \eqref{convo3} and noting that $|\Omega|/N\leq h^3$, we obtain
	\begin{equation}\label{lem1}
	|\partial_{x_i}\phi_{\alpha}|^2\leq\frac{C|\Omega|h^{-2}}{N}\|\partial_{x_i}\rho\|_{L^\infty(\Omega)}^2\leq Ch\|\partial_{x_i}\rho\|_{L^\infty(\Omega)}^2
	\end{equation}
\end{proof}
\noindent
We now prove the following theorem. The proof technique was borrowed from \cite{AOT}, where it was used to prove a similar statement for the two dimensional case.
\begin{theorem}
	Let $\mathcal{K}=\{-1,0,1\}^3$ and $\tilde{I}_h$ be as in \eqref{modef} and $\displaystyle K_{\rho}=\left(\sum_{i=1}^3\|\partial_{x_i}\rho\|_{L^\infty(\Omega)}^2\right)^{1/2}$. Then ,
	\begin{equation}\label{lem}
	|\tilde{I}_h(v)|_{H^1(\Omega)}^2 \leq ChK_{\rho}^2\sum_{\alpha\in\mathcal{J}}|\bar{v}_{\alpha}|^2~~ \forall v\in H^1(\Omega).
	\end{equation}
\end{theorem}
\begin{proof}
    We set $\epsilon=h/10.$ Hence, it follows immediately that $|\phi_{\alpha}\phi_{\beta}|=0$ for $\alpha-\beta\not\in \mathcal{K}$, where $\mathcal{K}=\{-1,0,1\}^3$.
	Differentiating \eqref{modef} and applying Cauchy-Schwarz inequality, we obtain
	\begin{equation}\label{aot}
	\begin{split}
	|\partial_{x_i}\tilde{I}_h(v)|^2&\leq \int_{\Omega}\sum_{\gamma\in\mathcal{K}}\sum_{\alpha\in\mathcal{J}}\left|\bar{v}_{\alpha}\partial_{x_i}\phi_{\alpha}(x)\right|\left|\bar{v}_{\alpha+\gamma}\partial_{x_i}\phi_{\alpha+\gamma}(x)\right|dx\\
	&\leq \int_{\Omega}\sum_{\gamma\in\mathcal{K}}\left(\sum_{\alpha\in\mathcal{J}}|\bar{v}_{\alpha}|^2\left|\partial_{x_i}\phi_{\alpha}(x)\right|^2\right)^{1/2}\left(\sum_{\alpha\in\mathcal{J}}|\bar{v}_{\alpha+\gamma}|^2\left|\partial_{x_i}\phi_{\alpha+\gamma}(x)\right|^2\right)^{1/2}dx\\
	&\leq  \int_{\Omega}\sum_{\gamma\in\mathcal{K}}\sum_{\alpha\in\mathcal{J}}|\bar{v}_{\alpha}|^2\left|\partial_{x_i}\phi_{\alpha}(x)\right|^2dx\\
	&\leq 27\int_{\Omega}\sum_{\alpha\in\mathcal{J}}|\bar{v}_{\alpha}|^2\left|\partial_{x_i}\phi_{\alpha}(x)\right|^2dx\\
	&\leq 27\sum_{\alpha\in\mathcal{J}}|\bar{v}_{\alpha}|^2|\partial_{x_i}\phi_{\alpha}|^2
	\end{split}
	\end{equation}
	Applying \eqref{lem1} to \eqref{aot}, we may write
	\begin{equation}\label{thm1}
	|\partial_{x_i}\tilde{I}_h(v)|^2 \leq\frac{C|\Omega|h^{-2}}{N}\|\partial_{x_i}\rho\|_{L^\infty(\Omega)}^2\sum_{\alpha\in\mathcal{J}}|\bar{v}_{\alpha}|^2\leq Ch\|\partial_{x_i}\rho\|_{L^\infty(\Omega)}^2\sum_{\alpha\in\mathcal{J}}|\bar{v}_{\alpha}|^2
	\end{equation}
	Now, summing over $i$, we obtain
	\begin{equation}\label{need1}
	|\tilde{I}_h(v)|_{H^1(\Omega)}^2\leq\frac{C|\Omega|h^{-2}}{N}K_{\rho}^2\sum_{\alpha\in\mathcal{J}}|\bar{v}_{\alpha}|^2 \leq ChK_{\rho}^2\sum_{\alpha\in\mathcal{J}}|\bar{v}_{\alpha}|^2~~ \forall v\in H^1(\Omega).
	\end{equation}
\end{proof}
\begin{corollary}
	Let $\tilde{I}_h$ be as in \eqref{modef} and $\displaystyle K_{\rho}=\left(\sum_{i=1}^3\|\partial_{x_i}\rho\|_{L^\infty(\Omega)}^2\right)^{1/2}$. Then
	\begin{equation}
	\|\tilde{I}_h(v)\| \leq CK_{\rho}\|v\|~~ \forall v\in H^1(\Omega).
	\end{equation}
\end{corollary}
\begin{proof}
	We first look at the term $\displaystyle\sum_{\alpha\in\mathcal{J}}|\bar{v}_{\alpha}|^2$ in \eqref{thm1}. Using the definition of $\bar{v}_{\alpha}$ given in \eqref{valpha}, Holder's inequality and Gagliardo-Nirenberg-Sobolev inequality, we obtain
	\begin{equation}\label{cor1}
	\begin{split}
	\sum_{\alpha\in\mathcal{J}}|\bar{v}_{\alpha}|^2 \leq\sum_{\alpha\in\mathcal{J}}\left(\left|\frac{1}{|Q_\alpha|}\int_{Q_{\alpha}}v(x)dx\right|\right)^2&\leq \frac{N^2}{|\Omega|^2}\sum_{\alpha\in\mathcal{J}}\|v\|^2_{L^1(Q_{\alpha})}\\
	&\leq \frac{N^2}{|\Omega|^2}\sum_{\alpha\in\mathcal{J}}|Q_\alpha|^{5/3}\|v\|^2_{L^6(Q_{\alpha})}\\
	&\leq \frac{CN^{1/3}}{|\Omega|^{1/3}}\sum_{\alpha\in\mathcal{J}}\|v\|^2_{H^1(Q_{\alpha})}\\
	&\leq \frac{CN^{1/3}}{|\Omega|^{1/3}}\sum_{\alpha\in\mathcal{J}}\sum_{i=1}^3\left\|\left(\partial_{x_i}v\right)\chi_{Q_{\alpha}}\right\|^2_{L^2(\Omega)}\\
	&\leq \frac{CN^{1/3}}{|\Omega|^{1/3}}\sum_{i=1}^3\left(\sum_{\alpha\in\mathcal{J}}\left(\partial_{x_i}v\right)\chi_{Q_{\alpha}},\sum_{\alpha\in\mathcal{J}}\left(\partial_{x_i}v\right)\chi_{Q_{\alpha}}\right)\\
	\end{split}
	\end{equation}
	Next, using the fact that for $\alpha,\beta\in\mathcal{J}$,
	\begin{equation}
	\chi_{Q_{\alpha}}\chi_{Q_{\beta}}=
	\begin{cases} 
	0 &\mbox{if } \alpha\neq\beta \\
	\chi_{Q_{\alpha}} & \mbox{if } \alpha=\beta ,
	\end{cases} 
	\end{equation}
	we can simplify the expression on the RHS of \eqref{cor1} to obtain
	\begin{equation}\label{h1}
	\begin{split}
	\sum_{\alpha\in\mathcal{J}}|\bar{v}_{\alpha}|^2&\leq \frac{CN^{1/3}}{|\Omega|^{1/3}}\sum_{i=1}^3\|\partial_{x_i}v\|^2_{L^2(\Omega)}\\
	&\leq \frac{CN^{1/3}}{|\Omega|^{1/3}}\|v\|^2
	\end{split}
	\end{equation}
	Applying \eqref{h1} to the the first inequality in \eqref{need1} and summing over $i$, we obtain
	\begin{equation}
	\|\tilde{I}_h(v)\| \leq CK_{\rho}\|v\|~~ \forall v\in H^1(\Omega).
	\end{equation}
\end{proof}
Next, we briefly try to see why $\tilde{I}_h$ is a type-\rom{1} interpolant.
We first look at a lemma, which is a modified version of a similar result in \cite{Evans}.
\begin{lemma}\label{p2}
		Let $U=\left\{(p,q,r)\in\mathbb{R}^3:|p-a|<0.5h, |q-b|<0.5h, |r-c|<0.5h, \text{ for }a,b,c\in \mathbb{R}\right\}$ be a cube of side length $h>0$ and center $m=(a,b,c)$ and $u\in W^{1,p}(U)$ . Assume $1\leq p\leq\infty$. Then there exists a $C$, depending on only on $p$, such that $$\|u-(u)_{U}\|_p\leq Ch\|\nabla u\|_p,$$ where $$\displaystyle(u)_U=\frac{1}{|U|}\int_{U}u(x)dx$$
		\begin{proof}
			From Theorem\ref{poincare},  it follows that the statement is true for a cube $V$ of side length one given by $V=\{(p,q,r)\in\mathbb{R}^3 : |p|<0.5, |q|<0.5, |r|<0.5\}$. For $y\in V$, we define $v\in W^{1,p}(V)$ by $$v(y)=u(hy+m).$$ From Theorem \ref{poincare} we obtain $$\|v-(v)_{V}\|_p\leq C\|\nabla v\|_p.$$
			Changing variables, we obtain the statement of the theorem.
		\end{proof}
	\end{lemma}
\begin{theorem}
    Let $\tilde{I}_h$ be the smoothed volume interpolant as defined in \eqref{modef}. Then $\tilde{I}_h$ is a type-\rom{1} interpolant.
\end{theorem}
\begin{proof}
    \begin{equation}\label{t1}
    \begin{split}
        |\tilde{I}_h(v)|^2&\leq\left|\sum_{\alpha\in\mathcal{J}}\bar{v}_\alpha\phi_\alpha\right|^2\\
    \end{split}
\end{equation}
Recall that $|\phi_{\alpha}\phi_{\beta}|=0$ for $\alpha-\beta\not\in \mathcal{K}$, where $\mathcal{K}=\{-1,0,1\}^3$. Repeating the arguments in \eqref{aot} and applying Young's inequality for convolutions, we obtain
\begin{equation}\label{i1}
\begin{split}
    |\tilde{I}_h(v)|^2&\leq 27\sum_{\alpha\in\mathcal{J}}|\bar{v}_{\alpha}|^2|\phi_{\alpha}|^2\\
    &\leq 27\sum_{\alpha\in\mathcal{J}}|\bar{v}_{\alpha}|^2\|\rho_{\epsilon}\|_{L^{1}(\Omega)}^2|\psi_{Q_\alpha}|^2
\end{split}   
\end{equation}
Observe that 
\begin{equation}\label{l1}
    \begin{split}
        \|\rho_{\epsilon}\|_{L^{1}(\Omega)}|^2&=\left(\int_{|x|\leq\epsilon}\left|\epsilon^{-3}\rho(x/\epsilon)\right|dx\right)^2\\
	&\leq C\|\rho\|_{L^\infty(\Omega)}^2
    \end{split}
\end{equation}
Applying \eqref{l1} and the middle inequality of \eqref{convo3} to \eqref{i1}, we obtain
\begin{equation}\label{the1}
    \begin{split}
        |\tilde{I}_h(v)|^2&\leq C\|\rho\|_{L^\infty(\Omega)}^2\frac{|\Omega|}{N}\sum_{\alpha\in\mathcal{J}}|\bar{v}_{\alpha}|^2\\
        &\leq C\|\rho\|_{L^\infty(\Omega)}^2\frac{|\Omega|}{N}\sum_{\alpha\in\mathcal{J}}\left(\left|\frac{1}{|Q_\alpha|}\int_{Q_{\alpha}}v(x)dx\right|\right)^2\\
        &\leq C\|\rho\|_{L^\infty(\Omega)}^2\frac{N}{|\Omega|}\sum_{\alpha\in\mathcal{J}}\|v\|^2_{L^1(Q_{\alpha})}\\
        &\leq C\|\rho\|_{L^\infty(\Omega)}^2\sum_{\alpha\in\mathcal{J}}\|v\|^2_{L^2(Q_{\alpha})}\\
        &\leq C\|\rho\|_{L^\infty(\Omega)}^2|v|^2
\end{split}
\end{equation}
We now look at $|v-\tilde{I}_h(v)|.$
\begin{equation}
\begin{split}
    |v-\tilde{I}_h(v)|^2\leq \left|\sum_{\alpha\in\mathcal{J}}(v-\bar{v}_\alpha)\phi_\alpha\right|^2
\end{split}
\end{equation}
After repeating  the arguments in \eqref{i1} and \eqref{l1}, we obtain
\begin{equation}
    \begin{split}
         |v-\tilde{I}_h(v)|^2&\leq C\|\rho\|_{L^\infty(\Omega)}^2\sum_{\alpha\in\mathcal{J}}|(v-\bar{v}_{\alpha})\phi_\alpha|^2\\
        &\leq C\|\rho\|_{L^\infty(\Omega)}^2\frac{|\Omega|}{N}\sum_{\alpha\in\mathcal{J}}|(v-\bar{v}_{\alpha})|^2\|\phi_\alpha\|_{L^{\infty}(\Omega)}^2\\
        &\leq C\|\rho\|_{L^\infty(\Omega)}^2\frac{N}{|\Omega|}\sum_{\alpha\in\mathcal{J}}|(v-\bar{v}_{\alpha})|^2\|\rho_{\epsilon}\|_{L^{\infty}(\Omega)}^2|\psi|_{L^1(\Omega)}^2\\
        &\leq Ch^{-6}\frac{|\Omega|^2}{N^2}\|\rho\|_{L^\infty(\Omega)}^4\sum_{\alpha\in\mathcal{J}}|(v-\bar{v}_{\alpha})|^2
    \end{split}
\end{equation}
Lastly, applying Lemma \ref{p2} and noting that $|\Omega|/N\leq h^3$, we obtain
\begin{equation}\label{the2}
\begin{split}
    |v-\tilde{I}_h(v)|^2&\leq Ch^2\|\rho\|_{L^\infty(\Omega)}^4\sum_{\alpha\in\mathcal{J}}\|v\|_{H^1(Q_\alpha)}^2\\
    &\leq Ch^2\|\rho\|_{L^\infty(\Omega)}^4\|v\|_{H^1(\Omega)}^2
\end{split}
\end{equation}
Therefore, from \eqref{the1} and \eqref{the2}, we see that $\tilde{I}_h$ is a type-\rom{1} interpolant.
\end{proof}
We now state and prove the following lemma due to E.S Titi (via private communication).
\begin{lemma}\label{titi}
	    Let $y:\mathbb{R}\to\mathbb{R}$ be a real valued function, with $y(t)\geq 0~\forall t\in \mathbb{R}$. If $y$ satisfies, for $0\leq s\leq t$ and a constant $\mu>0$, the inequality
	    \begin{equation}\label{1lemma}
	        y(t)+\mu\int_s^ty(\tau)d\tau\leq y(s),
	    \end{equation}
	    then
	    \begin{equation}\label{1lemma1}
	        y(t)\leq\mu e^{-\mu\left(t-\frac{1}{\mu}\right)}\int_{0}^{\frac{1}{\mu}}y(\tau)d\tau,~~\forall t\geq1/\mu.
	    \end{equation}
	\end{lemma}
	\begin{proof}
	    Let \begin{equation}
	        \displaystyle\psi(t)=\int_{t-\frac{1}{\mu}}^t y(\tau)d\tau.
	    \end{equation}
	    Choosing $s=\left(t-1/\mu\right)$,we can then rewrite \eqref{1lemma} as
	    \begin{equation}
	        \frac{\partial \psi(t)}{dt} +\mu\psi(t)\leq 0.
	    \end{equation}
	    Thus, from standard Gronwall over the interval $[\sigma,t]$ inequality, we obtain
	    \begin{equation}\label{11lemma2}
	        \psi(t)\leq\psi\left(\sigma\right)e^{-\mu\left(t-\sigma\right)}.
	    \end{equation}
	    Choosing $\sigma=1/\mu$, we obtain
	    \begin{equation}\label{1lemma2}
	        \psi(t)\leq\psi\left(1/\mu\right)e^{-\mu\left(t-\frac{1}{\mu}\right)}.
	    \end{equation}
	    Now, keeping only the first term on the LHS of \eqref{1lemma}, integrating with respect to $s$ over the interval $\left[t-\frac{1}{\mu},t\right]$ and applying \eqref{1lemma2}, we see that
	    \begin{equation}
	    \begin{split}\label{1lemma3}
	       \left(\frac{1}{\mu}\right) y(t)&\leq\int_{t-\frac{1}{\mu}}^t y(\tau)d\tau\\
	       &\leq \psi\left(1/\mu\right)e^{-\mu\left(t-\frac{1}{\mu}\right)}\\
	       &\leq e^{-\mu\left(t-\frac{1}{\mu}\right)}\int_0^{\frac{1}{\mu}}y(\tau)d\tau.
	    \end{split}
	    \end{equation}
	    Multiplying \eqref{1lemma3} by a factor of $\mu$, we obtain \eqref{1lemma1}.
	\end{proof}
	If, in \eqref{11lemma2}, we choose $\sigma=t/2$, we obtain the following Corollary
	\begin{corollary}\label{1lemmaaa1}
	      Let $y:\mathbb{R}\to\mathbb{R}$ be a real valued function, with $y(t)\geq 0~\forall t\in \mathbb{R}$. If $y$ satisfies, for $0\leq s\leq t$ and a constant $\mu>0$, the inequality
	    \begin{equation}
	        y(t)+\mu\int_s^ty(\tau)d\tau\leq y(s),
	    \end{equation}
	    then
	    \begin{equation}
	        y(t)\leq\mu e^{-\mu\left(t/2\right)}\int_{t/2-1/\mu}^{t/2}y(\tau)d\tau,~~\forall t\geq1/\mu.
	    \end{equation}
	\end{corollary}

\end{document}